\documentclass[%
 aip,
 amsmath,amssymb,
 reprint,%
]{revtex4-1}

\usepackage{graphicx}
\usepackage{dcolumn}
\usepackage{bm}

\usepackage[utf8]{inputenc}
\usepackage[T1]{fontenc}
\usepackage{mathptmx}
\usepackage{etoolbox}

\usepackage{placeins}
\usepackage{rotating}
\usepackage{svg}
\usepackage{amsmath}
\usepackage{float}
\usepackage{booktabs}
\usepackage{graphicx}

\usepackage{soul}
\usepackage{xcolor}

\makeatletter
\def\@email#1#2{%
 \endgroup
 \patchcmd{\titleblock@produce}
  {\frontmatter@RRAPformat}
  {\frontmatter@RRAPformat{\produce@RRAP{*#1\href{mailto:#2}{#2}}}\frontmatter@RRAPformat}
  {}{}
}%
\makeatother
\begin{document}

\preprint{AIP/123-QED}

\title[Evaluating Single-step Exponential Methods for Stiff ODEs]{Performance Evaluation of Single-step Explicit Exponential Integration Methods on Stiff Ordinary Differential Equations}
\author{Colby Fronk}
\affiliation{Department of Chemical Engineering; University of California, Santa Barbara; Santa Barbara, CA 93106, USA}
 \altaffiliation{Correspond to colbyfronk@ucsb.edu}
 
\author{Linda Petzold}%
\affiliation{Department of Mechanical Engineering and Computer Science; University of California, Santa Barbara; Santa Barbara, CA 93106, USA}%
\affiliation{Department of Computer Science; University of California, Santa Barbara; Santa Barbara, CA 93106, USA}
 \altaffiliation{Correspond to petzold@engineering.ucsb.edu}

\date{\today}

\begin{abstract}
Stiff systems of ordinary differential equations (ODEs) arise in a wide range of scientific and engineering disciplines and are traditionally solved using implicit integration methods due to their stability and efficiency.  However, these methods are computationally expensive, particularly for applications requiring repeated integration, such as parameter estimation, Bayesian inference, neural ODEs, physics-informed neural networks, and MeshGraphNets. Explicit exponential integration methods have been proposed as a potential alternative, leveraging the matrix exponential to address stiffness without requiring nonlinear solvers. This study evaluates several state-of-the-art explicit single-step exponential schemes against classical implicit methods on benchmark stiff ODE problems, analyzing their accuracy, stability, and scalability with step size.  Despite their initial appeal, our results reveal that explicit exponential methods significantly lag behind implicit schemes in accuracy and scalability for stiff ODEs.  The backward Euler method consistently outperformed higher-order exponential methods in accuracy at small step sizes, with none surpassing the accuracy of the first-order integrating factor Euler method. Exponential methods fail to improve upon first-order accuracy, revealing the integrating factor Euler method as the only reliable choice for repeated, inexpensive integration in applications such as neural ODEs and parameter estimation.  This study exposes the limitations of explicit exponential methods and calls for the development of improved algorithms.
\end{abstract}

\maketitle

\begin{quotation}
The numerical integration of stiff systems of ordinary differential equations (ODEs) is a common challenge in numerous scientific and engineering applications. Stiffness arises when a system exhibits rapid changes in some variables while others evolve slowly, as seen in chemical kinetics, electrical circuits, and population dynamics. Classical implicit integration methods, such as backward differentiation and implicit Runge-Kutta schemes, have long been the standard approach for dealing with stiff ODEs due to their ability to take larger timesteps while maintaining stability. However, these methods require solving nonlinear systems of equations at each timestep, which is both computationally expensive and time-consuming. In recent years, explicit exponential integration methods have been developed, aiming to circumvent these issues by leveraging matrix exponentials to account for the stiffness without the need for nonlinear solvers. Although these methods have shown promise, they have not yet seen widespread adoption, and there has been little direct comparison to classical implicit methods. This paper seeks to address this gap by evaluating a range of explicit single-step exponential integration methods. We assess their performance on several stiff ODE problems, comparing their stability and accuracy to that of traditional implicit methods.

\end{quotation}

\section{Introduction}

The challenge of integrating stiff ordinary differential equations (ODEs) originally stemmed from the need to effectively simulate processes such as chemical kinetics, where the equations exhibited significant stiffness because of varying timescales \cite{ascher1998computer}. To deal with this, numerical analysts have extensively researched and refined different approaches. Over time, implicit methods such as Backward Differentiation Formulas (BDF) and Radau schemes emerged as the favored techniques for stiff ODEs \cite{ascher1998computer}. These methods excelled in maintaining stability while allowing for larger time steps, outperforming explicit methods that required impractically small steps.  As integration methods became faster, the focus shifted to the challenge of efficiently solving the inverse problem, which involves identifying a model's parameter values through numerous repeated integrations with varying parameter values. Similarly, Bayesian inference \cite{chen2012monte, liang2011advanced, mcelreath2018statistical} involves sampling different parameter values in the ODE model, and each of these samples requires numerical integration to be carried out. The emergence of deep learning models that incorporate differential equations, such as neural ODEs \cite{chen2018neural, latent_ODEs, bayesianneuralode, stochastic_neural_ode, kidger2020neural, kidger2022neural, morrill2021neural, jia2019neural, chen2020learning, dagstuhl, doi:10.1063/5.0130803, fronk2024bayesian, fronk2024trainingstiffneuralordinary}, physics-informed neural networks \cite{owhadi2015bayesian, hiddenphysics, raissi2018numerical, raissi2017physics, raissi2017physics, osti_1595805, cuomo2022scientific, cai2021physics} (PINNs), and mesh-based frameworks like MeshGraphNets \cite{pfaff2021learningmeshbasedsimulationgraph}, has heightened the need for faster ODE solvers as well as the need for differentiable ODE solvers. Beyond training these deep learning approaches, the application of neural differential equations in control systems and robotics will require real-time inexpensive reliable numerical integration.

Although implicit schemes have been refined to reduce the computational cost of a single ODE integration, their speed is still not sufficient enough for tasks that require repeated integration. Explicit exponential integration methods \cite{trefethen2000spectral, boyd2001chebyshev, canuto1988spectral, ascher2008numerical, miranker2001numerical} have been proposed as a promising option for applications requiring repeated, low-cost integration. By leveraging matrix exponentials, these methods aim to handle stiffness without the computational burden of iterative solvers. While they are an attractive alternative, their practical performance in terms of accuracy and stability has not been thoroughly evaluated.  This study systematically benchmarks explicit exponential methods against classical implicit schemes to assess their accuracy, stability, and scalability with step size. The results provide critical insights into their feasibility for stiff ODEs and their potential application in fields such as neural ODEs and Bayesian inference. This work also highlights the urgent need for improved solvers that combine computational efficiency with robust accuracy to meet the demands of modern computational science and scientific machine learning.

Despite their initial promise, our results reveal that explicit exponential methods lag behind implicit schemes in accuracy and scalability for stiff ODEs.  None of these methods surpass the accuracy of the first-order integrating factor method, with even higher-order exponential schemes underperforming backward Euler at small step sizes.  This study identifies IF Euler as the only viable exponential integration method for repeated, inexpensive integration in applications such as neural ODEs and parameter estimation, while emphasizing the need for improved algorithms to address these limitations.

\FloatBarrier
\section{Methods}

This section opens with an introduction to stiff ODEs and their specific challenges. We then transition into an overview of implicit schemes, which are most commonly used but computationally expensive due to the need for solving nonlinear equations. We proceed by discussing a selection of single-step exponential integration schemes selected for their low cost and stability.

\begin{figure*}
    \centering
    \includegraphics[width=0.95\linewidth]{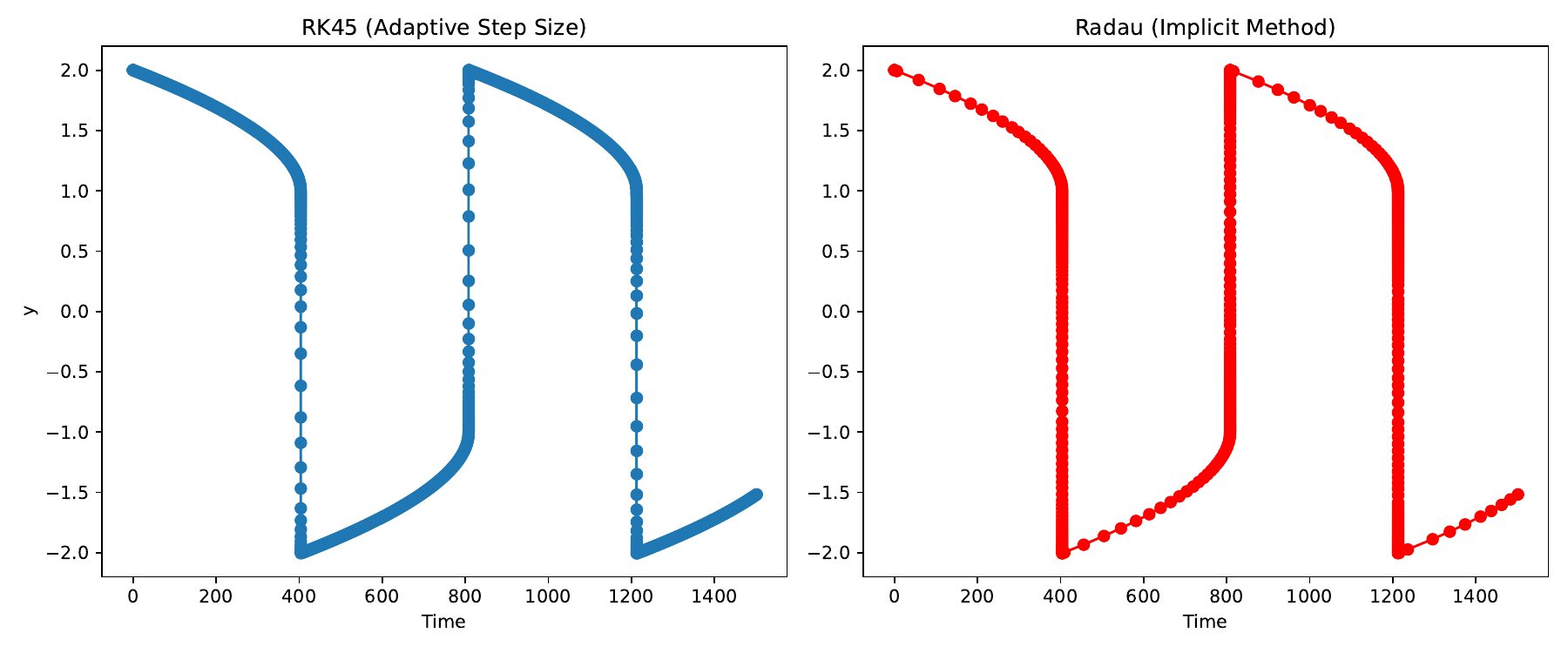}
    \caption{Comparison of the integration of the stiff van der Pol oscillator with $\mu=1000$ using two different integration schemes: explicit Runge-Kutta-Fehlberg, which is slow with 422,442 time points and 2,956,574 function evaluations, and implicit Radau IIA 5th order, which is faster with only 857 time points and 7,123 function evaluations. }
    \label{fig:stiffness_explained}
\end{figure*}

\subsection{Stiff ODEs}

Stiffness occurs when there is a large disparity in time scales, often observed when the Jacobian matrix's eigenvalues cover a wide range of magnitudes. When classical explicit integration methods such as the explicit Runge-Kutta method are used to integrate these stiff systems, a large number of time steps are taken to maintain numerical stability rather than accuracy requirements. In contrast, implicit schemes such as the Radau \cite{axelsson1969class, ehle1969pade, hairer1999stiff} methods are able to take larger time steps and maintain stability. We visually show this in Figure \ref{fig:stiffness_explained}. Explicit Runge-Kutta-Fehlberg integrates the stiff Van der Pol Oscillator with $\mu=1000$ using 422,442 time points and 2,956,574 function evaluations, while the implicit Radau IIA 5th order method is much faster with only 857 time points and 7,123 function evaluations. This example illustrates the core issue of stiffness in differential equations.  Classical explicit methods require very small time steps to achieve stability for stiff problems, greatly increasing computational cost.

\subsection{Implicit Methods are Expensive}

The integration of stiff ODEs is conventionally carried out using implicit schemes. Unlike explicit methods, which compute the next step based solely on the current state, implicit methods require solving nonlinear systems at each time step. While this demands greater computational resources, it offers significantly improved stability, enabling larger time steps without sacrificing accuracy. This makes them well-suited for stiff systems, though they are less efficient for non-stiff ODEs where stability is not a significant concern.

Some common implicit methods for solving stiff ODEs are the backward Euler method, the trapezoid method, and Radau IIA methods. The backward Euler method is the simplest implicit integration scheme:  

\begin{equation}
y_{n+1} = y_n + h f(t_{n+1}, y_{n+1}),
\end{equation}

\noindent Backward Euler is A-stable, providing strong damping of oscillations and stability for stiff systems. However, it can overdampen solutions, making them too smooth. The trapezoid method improves upon backward Euler by averaging the derivatives at the current and next points: 

\begin{equation}
y_{n+1} = y_n + \frac{h}{2} \left( f(t_n, y_n) + f(t_{n+1}, y_{n+1}) \right).
\end{equation}

\noindent Trapezoid method is second-order accurate and A-stable but not L-stable, meaning it may allow some stiff components to oscillate rather than decay smoothly. As a highly stable subclass of implicit Runge-Kutta methods, Radau IIA methods \cite{axelsson1969class, ehle1969pade, hairer1999stiff} have the following form:

\begin{align}
Y_i =& y_n + h \sum_{j=1}^{s} a_{ij} f(t_n + c_j h, Y_j), \quad i = 1, \ldots, s, \\
y_{n+1} =& y_n + h \sum_{j=1}^{s} b_{j} f(t_n + c_j h, Y_j)
\end{align}

\noindent where \(Y_i\) are the stage values, \(a_{ij}\) are the coefficients from the Runge-Kutta method's Butcher tableau, and \(c_j\) are the nodes the solution is found at. The Butcher tableau for Radau3 and Radau5 can be found in Ref.~\onlinecite{hairer1999stiff}. The Radau3 and Radau5 methods are both A-stable and L-stable, which ensures rapid decay of transient components without oscillations, even with large time steps. 

For all of these implicit methods discussed, the solution at the next step, $y_{n+1}$ is found by solving a multivariate nonlinear system of equations using iterative techniques such as Newton's method. This is why implicit schemes are more computationally costly per integration step compared to explicit schemes. Although the Radau IIA methods are higher-order methods, they require solving for \(Y_i\) a nonlinear solution involving a matrix-valued function, which introduces additional numerical challenges.

\subsection{Explicit Exponential Integration Methods}

Exponential integration methods \cite{trefethen2000spectral, boyd2001chebyshev, canuto1988spectral, ascher2008numerical, miranker2001numerical} approach the solution of an ODE by separating it into a linear component and the remaining nonlinear terms:

\begin{equation}
    \frac{dy}{dt} = f(t,y(t)) = L y + N(t,y)
    \label{eqn:linearizedODE}
\end{equation}

\noindent where the linear term, $L$, is the Jacobian evaluated at our initial condition $y_0$:

\begin{equation}
    L = \frac{df}{dy}(y_0)
\end{equation}

\noindent and the nonlinear term, $N$, represents the remaining components:

\begin{equation}
    N(t,y(t)) = f(t,y(t)) - Ly.
    \label{eqn:nonlinear}
\end{equation}

\noindent This decomposition allows for the exact integration of the linear component using the matrix exponential, which is particularly beneficial when the linear term represents the stiffest part of the dynamics, and the nonlinear part changes more slowly. By handling the stiff linear dynamics exactly and using an explicit method for the slower nonlinear part, we can better manage stiffness. However, this reformulation still leaves the question of how to effectively integrate the new ODE. The following sections discuss different integration techniques that leverage this exponential integration framework. 

\subsection*{Integrating Factor Methods}

Integrating factor methods \cite{lawson1967generalized, boyd2001chebyshev, canuto1988spectral, fornberg1999fast, trefethen2000spectral} use a change of variables

\begin{equation}
    w(t) = e^{-Lt} y(t),
\end{equation}

\noindent which, when differentiated with respect to $t$, results in:

\begin{equation}
    \label{eqn:IF_eqn1}
    \frac{dw(t)}{dt} = e^{-Lt} \left (\frac{dy(t)}{dt} - Ly(t) \right ).
\end{equation}

\noindent By substituting Eqn. \ref{eqn:nonlinear} into Eqn. \ref{eqn:IF_eqn1}, we obtain:

\begin{equation}
    \frac{dw(t)}{dt} = e^{-Lt} N(t,y(t)) = e^{-Lt} N(t,e^{Lt}y(t)).
    \label{eqn:integratingfactor}
\end{equation}

\noindent Integrating factor methods are straightforward to derive because they allow the use of any standard numerical integration scheme, such as multi-step or Runge-Kutta methods, to numerically solve Eqn \ref{eqn:integratingfactor}. After applying the chosen integration method, the solution is converted back to the original variable $y$. Despite their simplicity, these methods are often discouraged in the literature. They tend to be inaccurate, particularly when the nonlinear term changes slowly, and they can fail to compute fixed points correctly, resulting in greater errors compared to other exponential integration techniques.

\subsubsection{Integrating Factor Euler Method (IF Euler)}

Applying the Forward Euler method on Eqn \ref{eqn:integratingfactor} results in the A-stable first-order Integrating Factor Euler scheme \cite{lawson1967generalized, boyd2001chebyshev, canuto1988spectral, fornberg1999fast, trefethen2000spectral, cox2002exponential}, which has local truncation error $\frac{h^2}{2}LN$:

\begin{equation}
    y_{n+1} = e^{Lh} (  y_n + h N_n  ).
    \label{eqn:IFEuler}
\end{equation}

\subsubsection{Second-Order Integrating Factor Runge-Kutta Method (IF2RK)}

Applying the explicit second-order Runge-Kutta method to Eqn \ref{eqn:integratingfactor} gives us the second-order Integrating Factor Runge-Kutta scheme \cite{lawson1967generalized, boyd2001chebyshev, canuto1988spectral, fornberg1999fast, trefethen2000spectral, cox2002exponential}, which has local truncation error $\frac{5h^3}{12}L^2N$:

\begin{equation}
    y_{n+1} = e^{Lh} y_n + \frac{h}{2} \left[ e^{Lh} N_n + N(t_n+h, e^{Lh} (  y_n + h N_n  )) \right].
\end{equation}

\noindent The local truncation error of this method shows why the literature discourages the use of integrating factor exponential integration methods. Since the stiff components are in the linear term $L$, the $L^2$ term rapidly grows as the ODE gets stiffer.

\subsection*{Exponential Time Differencing (ETD) Schemes}

In exponential time differencing methods \cite{beylkin1998new, certaine1960solution, ascher2008numerical, miranker2001numerical}, we reformulate Eqn \ref{eqn:linearizedODE} as follows:

\begin{equation}
    y(t_{n+1}) = e^{Lh} y(t_n) + \int_0^h e^{-L \tau} N(t_n+\tau, y(t_n+\tau))  d\tau. 
    \label{eqn:ETD_generalform}
\end{equation}

\noindent The linear component is integrated exactly using the matrix exponential, while the integral involving the nonlinear part is approximated. This formulation allows for exact computation of the fixed points of the ODE, unlike the integrating factor methods, and offers improved accuracy. Due to these advantages, exponential time differencing methods are generally preferred in the literature over integrating factor techniques. However, deriving schemes based on this approach is considerably more complex and challenging.

\subsubsection{First-order ETD1 Scheme}
    The simplest approximation comes from approximating the nonlinear term to be a constant $N \approx N_n$ to give us the following first-order ETD1 scheme \cite{cox2002exponential} with local truncation error of $\frac{h^2}{2} \dot{N}$:
\begin{equation}
    y_{n+1} = e^{Lh} y_n + \left(e^{Lh} - I \right) \frac{N_n}{L}. 
\end{equation}

\subsubsection{Second-order Runge-Kutta ETD Scheme (ETD2RK)}
A second-order Runge-Kutta ETD scheme \cite{cox2002exponential} with local truncation error of $\frac{h^3}{12} \ddot{N}$ can be constructed as follows:
\begin{align}
    a_n =& e^{Lh} y_n + \left(e^{Lh} - I \right) \frac{N_n}{L} \\ \notag
    y_{n+1} =& a_n + \left(e^{Lh} - I - hL \right) \frac{N(t_n+h, a_n)-N_n}{hL^2}
\end{align}

\subsubsection{Classical Fourth-order Runge-Kutta ETD Scheme (ETD4RK)}

A direct extension of the classical fourth-order Runge-Kutta scheme on the exponential differencing equation, only gives a third-order scheme \cite{cox2002exponential}:

\begin{align}
    a_n =& e^{ L h / 2 } y_n + \left( e^{Lh/2} - I \right) \frac{N_n}{L}, \notag \\
    b_n =& e^{ L h / 2 } y_n + \left( e^{Lh/2} - I \right) \frac{N(t_n + h/2, a_n )}{L}, \notag \\
    c_n =& e^{ L h / 2 } a_n + \left( e^{Lh/2} - I \right) \frac{ 2 N(t_n + h/2, b_n ) - N_n}{L}, \notag \\
    T_1 =& \left[ -4 - Lh + e^{Lh} \left( 4 - 3 L h + (L h)^2 \right)  \right] N_n, \notag \\
    T_2 =& 2 \left[ 2 + L h + e^{Lh} \left( -2 + L h \right) \right] \left[ N(t_n+h/2, a_n ) \right. \notag \\
        & \left. + N( t_n + h / 2, b_n ) \right], \notag \\
    T_3 =& \left[ -4 - 3L h - (Lh)^2 + e^{Lh} \left(4 - Lh \right) \right] N(t_n + h, c_n ), \notag \\
    y_{n+1} =& e^{L h} y_n + \frac{T_1 + T_2 + T_3}{h^2L^3}.
\end{align}

Although this ETD4RK method is formally fourth-order, Ref.~\onlinecite{hochbruck2010exponential} demonstrated that it does not meet certain stiff-order conditions, effectively reducing its accuracy to third-order.

\subsubsection{Pseudo-Steady-State Approximation (PSSA) scheme (RKMK2e)}

The Pseudo-Steady-State Approximation (PSSA) scheme \cite{verwer1994evaluation, maset2009unconditional} is given by:

\begin{align}
    a_n &= e^{L_n h} y_0 + \frac{e^{L_n h} - I}{L} N_n(t_0, y_0), \\ 
    y_{n+1} &= e^{L_n h} y_0 + \frac{1}{2} \frac{e^{L_n h} - I}{L} \left( N_n(t_0 + h, a_n) + N_n(t_0, y_0) \right). \notag
\end{align}

\noindent This scheme belongs to the class of exponential Runge-Kutta schemes.

\subsubsection{Second-order ETD-real Distinct Poles Scheme (ETD-RDP)}

Computing the matrix exponential can be computationally expensive. To avoid this, we can use a second-order rational approximation for the matrix exponential, allowing us to reformulate the ETD scheme without explicitly calculating the matrix exponential \cite{iyiola2018exponential}:

\begin{equation}
    r(z)=\frac{1-\frac{5}{12}z}{(1+\frac{1}{3}z)(1+\frac{1}{4}z)} = \frac{9}{1+\frac{1}{3}z} - \frac{8}{1+\frac{1}{4}z} \approx e^{-z}.
\end{equation}

\noindent With this approximation the following second-order ETD-real distinct poles scheme \cite{iyiola2018exponential} can be obtained:

\begin{align}
    y^* &= \left( I + Lh  \right)^{-1} \left( y_n + hN_n  \right)\\ \notag
    y_{n+1} &= \left( I+\frac{1}{3}Lh  \right)^{-1} \left( 9y_n + 2hN_n +hN(t_n+h,y^*)\right) \\\notag  &+ \left( I+\frac{1}{4}Lh  \right)^{-1} \left( -8y_n - \frac{3h}{2}N_n - \frac{h}{2} N(t_n+h,y^*)  \right)
\end{align}

\subsection*{Exponential Propagation Iterative (EPI) Schemes}

Exponential propagation iterative (EPI) schemes \cite{tokman2006efficient} are constructed by expanding the right hand side of the ode $f(t,y)$ in the following way: 

\begin{equation}
    f(t,y(t)) = f_n + L \left( y(t) - y_n \right) + R(y(t))
\end{equation}

\begin{equation}
    R(y(t)) = f(t,y(t)) - f(y_n) - L \left(y(t) - y_n \right)
\end{equation}

\noindent We can use an integrating factor $e^{-Lt}$ to obtain

\begin{equation}
    \frac{d}{dt}\left( e^{-Lt} y(t)\right) = e^{-Lt} \left( f_n - L y_n \right) + e^{-Lt} R(y(t))
\end{equation}

\noindent and integrate this equation over the time interval to yield the integral form:

\begin{equation}
    y(t_n+h) = y_n + \left( e^{Lh} - I \right) L^{-1} f_n + \int_{t_n}^{t_n+h} e^{L (t_n+h-t)} R(y(t))  dt.
    \label{eqn:EPI_generalform}
\end{equation}

\noindent Equation \ref{eqn:EPI_generalform} is the starting point for constructing various exponential propagation iterative (EPI) schemes.  EPI schemes are constructed in a similar way to ETD schemes.  The integral is approximated using various quadrature rules for the integral.  For the case when the stiffness of the problem comes exclusively from the linear terms of $f(t,y)$, Equation \ref{eqn:EPI_generalform} reduces to Equation \ref{eqn:ETD_generalform} and we have the ETD class of integration schemes.  

\subsubsection{Second-order exponential propagation iterative scheme (EPI2)}

The simplest EPI scheme is the second-order exponential propagation iterative scheme \cite{tokman2006efficient} (EPI2), which is given by:

\begin{equation}
    y_{n+1} = y_n + \frac{e^{Lh}-I}{L} f_n.
\end{equation}

\subsubsection{Third-order Runge-Kutta type exponential propagation iterative scheme (EPIRK3)}

Using Runge-Kutta as our quadrature rule gives the following third-order Runge-Kutta type exponential propagation iterative scheme \cite{tokman2006efficient} (EPIRK3):

\begin{align}
    r_1 &= y_n + 2 \frac{e^{L\frac{h}{2}}-I}{L} f_n \\ \notag
    y_{n+1} &= y_n + \frac{e^{Lh}-I}{L} f_n + \frac{1}{3} \frac{e^{Lh}-(I+Lh)}{L^2h} R(r_1).
\end{align}

\subsection*{Generalized integrating factor exponential methods}

Generalized integrating factor methods \cite{krogstad2005generalized, li2023low, maday1990operator} extend the traditional integrating factor approach previously discussed. These methods simplify the ODE by finding a modified vector field to approximate the original vector field near a specific point, allowing the ODE to be solved exactly while capturing the essential features of the system. Additionally, interpolating polynomial approximations are used to approximate the nonlinear terms. Unlike standard integrating factors, generalized integrating factors preserve the fixed points of the system. 

\subsubsection{ETD1/RK4}

ETD1/RK4 is a generalized integrating factor exponential method that combines first-order exponential time differencing (ETD) with the RK4 method \cite{krogstad2005generalized, li2023low, maday1990operator}:

\begin{align}
    \phi_1(z) &= \frac{e^z-I}{z}\\ \notag
    \phi_{\tau, \Tilde{F}}(v) &= e^{\tau L}v + \tau\phi_1(\tau L)N_n 
\end{align}

\begin{align*}
    a =& \phi_{h/2,\Tilde{F}}(y_n) & N_a =& N(t_n+h/2, a) \\
    b =& \phi_{h,\Tilde{F}}(y_n) &  \\
    c =& a+\frac{h}{2} \left( N_a - N_n \right) & N_c =& N(t_n+h/2, c) \\
    d =& b+he^{\frac{h}{2}L} \left( N_c -N_n \right) & N_d =& N(t_n+h,d)
\end{align*}
\begin{equation*}
    y_{n+1} = b + \frac{h}{3} e^{\frac{h}{2}L} \left( N_a+N_c-2N_n \right) + \frac{h}{6} \left( N_d - N_n \right)
\end{equation*}

\noindent The use of higher-order interpolating polynomials creates higher-order schemes such as ETD2/RK4 and ETD3/RK4.  However, there is a considerable decrease in stability as the order of the method increases. We focus on ETD1/RK4 in our study because it provides the highest level of stability.

\subsection*{Strong stability preserving Runge-Kutta methods}

Strong stability preserving Runge-Kutta methods \cite{isherwood2019strong, carrillo2003weno, shu1988total, shu1988efficient, ferracina2004stepsize, ferracina2005extension, ferracina2008strong, gottlieb2011strong, gottlieb1998total, gottlieb2001strong, higueras2004strong, higueras2005representations, hundsdorfer2003monotonicity, ketcheson2008highly, ketcheson2009optimal, ketcheson2009computation, ketcheson2011step, kubatko2014optimal} were created by the partial differential equation (PDE) community for greater stability in solving hyperbolic PDEs. Strong stability preserving (SSP) time discretizations \cite{shu1988total, shu1988efficient} were developed to improve the stability of nonlinear, non-inner-product spatial discretizations when used with the forward Euler method. These SSP methods were further extended to higher-order time discretizations.

\subsubsection{Explicit SSP Runge-Kutta method (eSSPRK)}

Explicit SSP Runge-Kutta method \cite{isherwood2019strong} is given by:

\begin{align}
    u_1 &= e^{L_n h } \left( y_0 +  h N_n(t_0, y_0) \right), \notag \\
    u_2 &= \frac{3}{4} e^{L_n h / 2} y_0 + \frac{1}{4} e^{-L_n h / 2}  \left( u_1 + h N_n(t_0, u_1) \right) , \notag \\
    y_{n+1} &= \frac{1}{3} e^{L_n h} y_0 + \frac{2}{3} e^{L_n h / 2}  \left( u_2 + h N_n(t_0, u_2) \right) .  \\
\end{align}

\subsubsection{Explicit SSP Runge-Kutta method with nondecreasing abscissas (eSSPRK\textsuperscript{+})}

Explicit SSP Runge-Kutta method with nondecreasing abscissas \cite{isherwood2019strong} is given by:

\begin{align}
    u_1 &= e^{2 L_n h / 3} \left( \frac{1}{2} y_0 + \frac{1}{2} \left( y_0 + \frac{4}{3} h N_n(t_0, y_0) \right) \right), \notag \\
    u_2 &= e^{2 L_n h / 3} \left( \frac{2}{3} y_0 + \frac{1}{3} \left( u_1 + \frac{4}{3} h N_n(t_0, u_1) \right) \right), \notag \\
    y_{n+1} &= e^{L_n h} \left( \frac{59}{128} y_0 + \frac{15}{128} \left( y_0 + \frac{4}{3} h N_n(t_0, y_0) \right) \right) \notag \\
    &\quad + \frac{27}{64} e^{L_n h / 3} \left( u_2 + \frac{4}{3} h N_n(t_0, u_2) \right).
\end{align}

\noindent This method is designed to be more stable than the explicit SSP Runge-Kutta method.

\subsection*{Computing the Matrix Exponential}

Computing the matrix exponential \cite{del2003survey, RUIZ2016370, ARIOLI1996111}, \( e^A \), for these exponential integration schemes can be computationally expensive. For dense matrices, a common approach is the scaling and squaring method with a Pade approximation, which involves scaling down the matrix, applying a rational approximation, and squaring the result. This method has a complexity of \( \mathcal{O}(n^3) \), making it impractical for large matrices. If \( A \) is diagonalizable, \( e^A \) can be computed as \( e^A = V e^\Lambda V^{-1} \), where \( \Lambda \) is diagonal, but this still has a complexity of \( \mathcal{O}(n^3) \) due to costly matrix multiplications and inversions. For large sparse matrices, Krylov subspace methods, such as Arnoldi iteration, can reduce the complexity to \( \mathcal{O}(n^2) \). However, these iterative methods may require numerous iterations to reach the desired accuracy, which can increase computational costs and reduce stability. As the dimension of the ODE system grows, the cost of computing the matrix exponential, which scales as \( \mathcal{O}(n^2) \) in the best case scenario for iterative methods, can become prohibitively expensive.

\FloatBarrier
\clearpage
\section{Results}

The models explored in this study are the Stiff Van der Pol Oscillator \cite{VanderPolModel}, the HIRES model \cite{Schaefer1975HIR}, and the Robertson model \cite{robertson1966reaction}. The two-dimensional Van der Pol Oscillator has a tunable stiffness parameter $\mu$, which we set to a value of 1000 to make the system have steep oscillatory transients followed by stiff regions with little change. The eight-dimensional HIRES model tests the performance of the integration schemes on a moderately stiff and moderately sized chemical reaction network. The Robertson model, known for high stiffness and large coefficient values, introduces the additional complexity of logarithmic time spacing to challenge integration over extremely large time steps, reaching orders of $10^7$.

To establish testing intervals for each model, we first solved the ODEs over time using the SciPy \cite{2020SciPy-NMeth} Radau solver at its default tolerance, achieving time points uniformly spaced with errors of a similar magnitude. To develop test configurations suited for a range of time step sizes, we either repeatedly halved the dataset by removing every other time point or repeatedly doubled it by adding midpoints between each consecutive pair. After obtaining the desired step size, we grouped the data into pairs of data points, each pair starting with the initial point $t_0$ and ending with the next point $t_1$.  For each data pair, we applied each single-step integration method to predict the solution at $t_1$ given the known value at $t_0$. To analyze the accuracy of each of the methods, we plotted the absolute error of each of the predictions vs time. In addition to the previously mentioned exponential integration schemes, we also tested the explicit RK4 scheme which is known to be unstable for stiff ODEs.  For comparison, we also tested the standard single-step implicit schemes of backward Euler, trapezoid method, Radau3, and Radau5.

\subsection{Example 1: Stiff Van der Pol Model}

The stiff Van der Pol oscillator \cite{VanderPolModel} is given by the following system of ODEs:

\begin{equation}
\label{eqn:example3}
\begin{aligned}
    &\frac{dx}{dt} = y,   \\
    &\frac{dy}{dt} = 1000 y - 1000 x^2 y - x,  \\
    &x(0) = 1, \quad y(0)=0, \quad t \in [0,1300].
\end{aligned}
\end{equation}

\begin{figure*}
    \centering
    \includegraphics[width=0.9\linewidth]{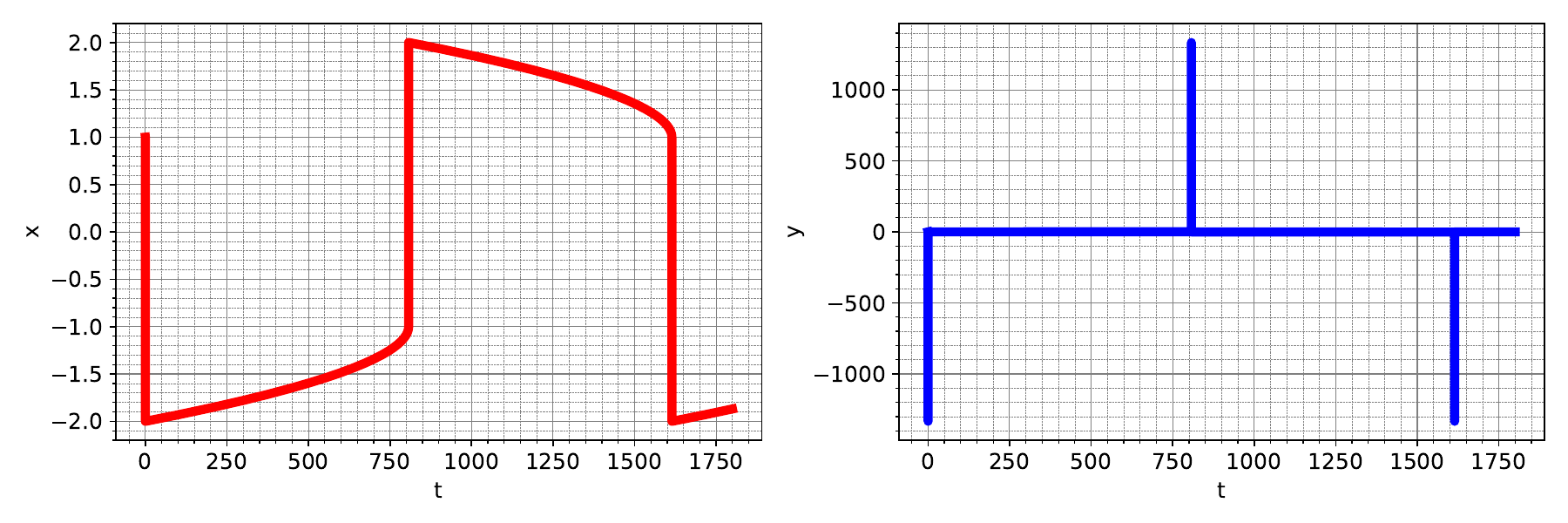}
    \caption{The stiff Van der Pol model. We have set $\mu=1000$, creating an oscillatory system with high stiffness. The model has sharp jumps in y, followed by stiff areas where the solution evolves more slowly.}
    \label{fig:vanderpol_data}
\end{figure*}

\noindent Within the field of dynamical systems, the Van der Pol oscillator is frequently cited as a classic example of a system exhibiting nonlinear damping and relaxation oscillations \cite{VanderPolModel}. Examples of its use include modeling neuronal action potentials, seismic fault movements, and vocal cord vibrations during speech \cite{FITZHUGH1961445, 4066548, doi:10.1142/S0218127499001620, doi:10.1121/1.4798467}. The stiffness in this model is dictated by the parameter $\mu$, which can be adjusted, making it a useful tool for evaluating the performance of stiff neural ODE solvers. In the equation above, we have set $\mu=1000$, creating an oscillatory system with high stiffness. As shown in Figure \ref{fig:vanderpol_data}, the model exhibits sharp jumps in y followed by stiff regions where the solution changes slowly. 

\begin{turnpage}
    \begin{figure*}
    \centering
    \includegraphics[width=0.95\linewidth]{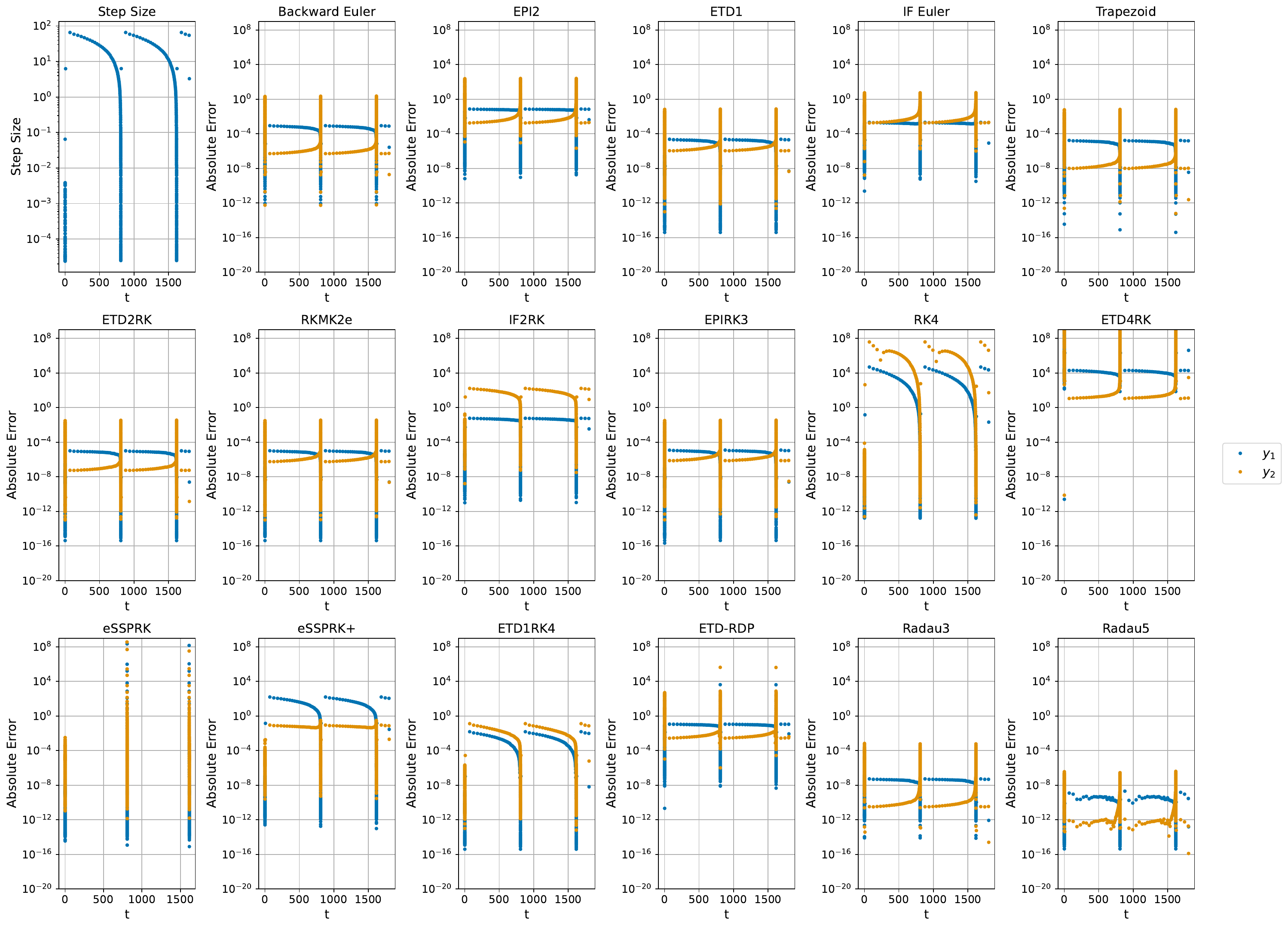}
    \caption{Absolute error for various integration schemes on the stiff Van der Pol model. The step size is also shown for reference. The data shown correspond to $n=1555$ data points across the integration time interval.}
    \label{fig:1555_VanderPol_plot}
\end{figure*}
\end{turnpage}

\begin{turnpage}
    \begin{figure*}
    \centering
    \includegraphics[width=0.95\linewidth]{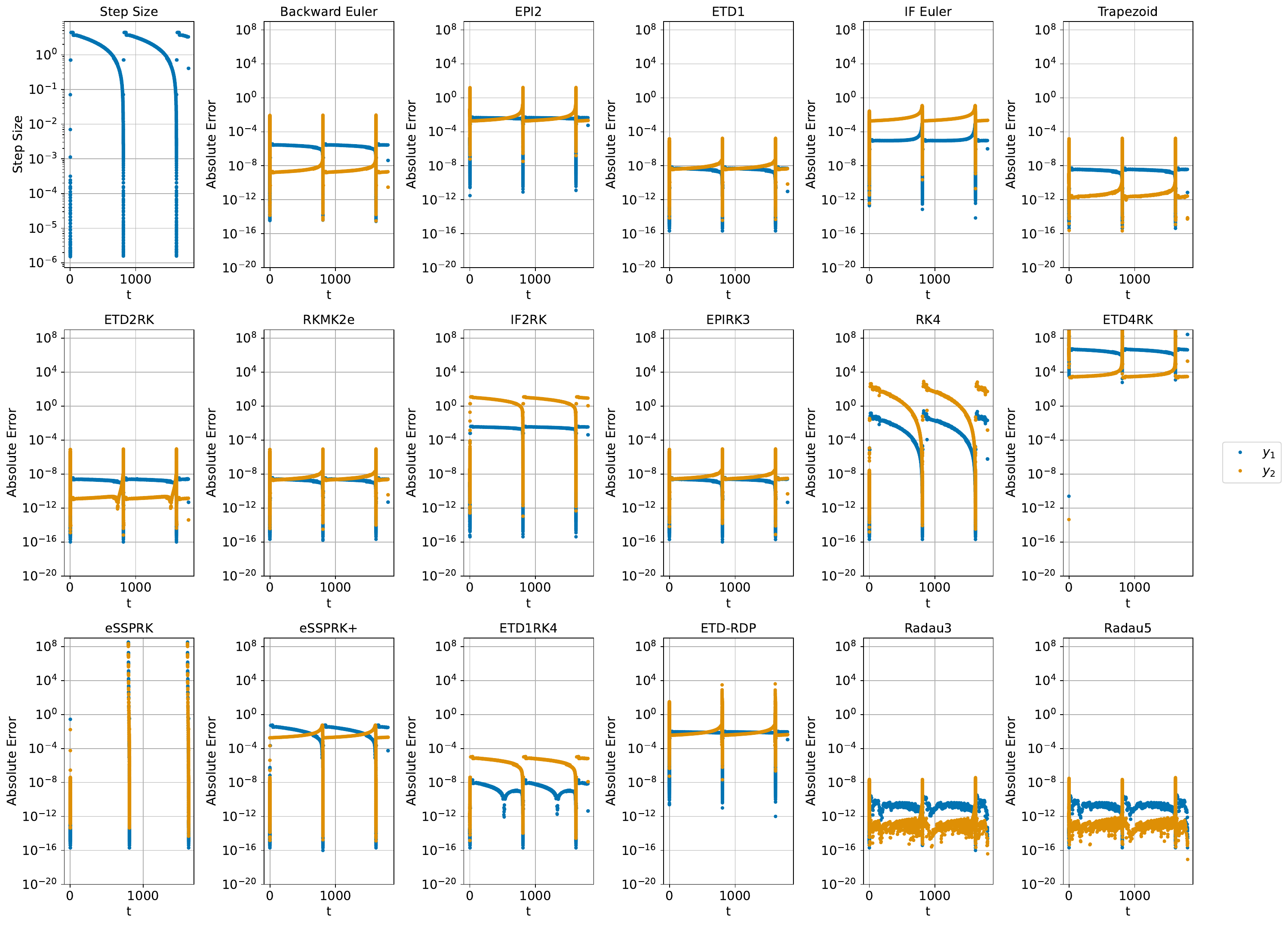}
    \caption{Absolute error for various integration schemes on the stiff Van der Pol model. The step size is also shown for reference. The data shown correspond to $n=24849$ data points across the integration time interval.}
    \label{fig:24849_VanderPol_plot}
\end{figure*}
\end{turnpage}

Figures \ref{fig:1555_VanderPol_plot} and \ref{fig:24849_VanderPol_plot} show the absolute error for the various integration schemes corresponding to 1555 and 24849 data points, respectively. For larger time steps ($n=1555$), the explicit RK4 integrator becomes unstable and produces significant errors, demonstrating the need for more robust methods that can address stiffness. The methods ETD4RK and eSSPRK perform worse than RK4. The implicit schemes backward Euler, Trapezoid, Radau3, and Radau5 are also shown for comparision purposes. Schemes EPI2, IF2RK, eSSPRK+, ETDIRK4, and ETD-RDP perform worse than backward Euler. eSSPRK+ appears to be more stable than eSSPRK, which is not surprising because it was designed to improve the stability of eSSPRK. IF Euler, ETD2RK, and RKMK2e all have about the same accuracy as backward Euler. Surprisingly, none of the explicit exponential methods come close to the accuracy of the Radau methods for large time steps. 

For smaller time steps ($n=24849$), RK4, ETD4RK, and eSSPRK are still unstable as shown by their large errors.  ETD4RK and eSSPRK still perform considerably worse than RK4. eSSPRK+ is more stable than eSSPRK as well.  EPI2, eSSPRK+, IF2RK, and ETD-RDP show worse performance than backward Euler. IF Euler method and backward Euler have comparable accuracy. The exponential methods with the highest accuracy are ETD1, ETD2RK, RKMK2e, EPIRK3, and ETDIRK4 with performance about the same as Trapezoid method. The exponential methods explored for the stiff Van der Pol model seem to saturate at an accuracy comparable to the Trapezoid method. Even at smaller time steps, none of the exponential methods have an accuracy close to the shown Radau methods.

\clearpage
\FloatBarrier
\subsection{Example 2: HIRES Model}

The “High Irradiance RESponse” (HIRES) model \cite{Schaefer1975HIR} is given by the following set of equations:

\begin{equation}
\label{eqn:hires_model}
\begin{aligned}
    &\frac{dy_1}{dt} = -1.71y_1 + 0.43y_2 + 8.32y_3 + 0.0007, \\
    &\frac{dy_2}{dt} = 1.71y_1 - 8.75y_2, \\
    &\frac{dy_3}{dt} = -10.03y_3 + 0.43y_4 + 0.035y_5, \\
    &\frac{dy_4}{dt} = 8.32y_2 + 1.71y_3 - 1.12y_4, \\
    &\frac{dy_5}{dt} = -1.745y_5 + 0.43y_6 + 0.43y_7, \\
    &\frac{dy_6}{dt} = -280y_6y_8 + 0.69y_4 + 1.71y_5 - 0.43y_6 + 0.69y_7, \\
    &\frac{dy_7}{dt} = 280y_6y_8 - 1.81y_7, \\
    &\frac{dy_8}{dt} = -280y_6y_8 + 1.81y_7, \\
    &y_1(0) = 1, \quad y_2(0)=0, \quad y_3(0)=0, \quad y_4(0)=0, \\
    &y_5(0) = 1, \quad y_6(0)=0, \quad y_7(0)=0, \quad y_8(0)=0.0057, \\
    &t \in [0, 321.8122].
\end{aligned}
\end{equation}

\noindent The HIRES model describes the process by which organisms adjust their growth and development in response to environmental light conditions \cite{Schaefer1975HIR}. The eight variables in the model represent concentrations of biochemical species. 

\begin{turnpage}
    \begin{figure*}
    \centering
    \includegraphics[width=0.95\linewidth]{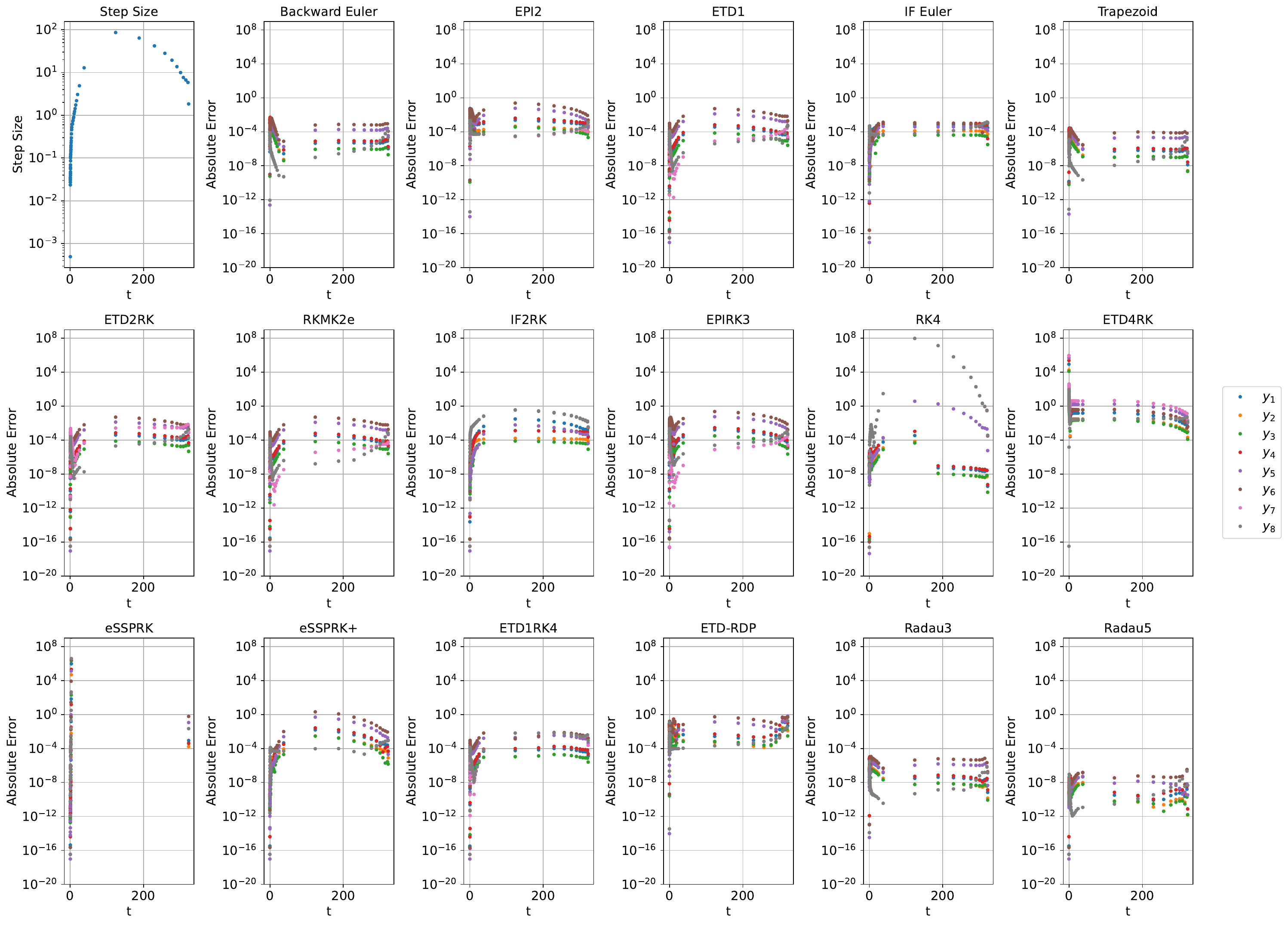}
    \caption{Absolute error for various integration schemes on the HIRES model. The step size is also shown for reference. The data shown correspond to $n=56$ data points across the integration time interval.}
    \label{fig:56_HIRES_plot}
\end{figure*}
\end{turnpage}

\begin{turnpage}
    \begin{figure*}
    \centering
    \includegraphics[width=0.95\linewidth]{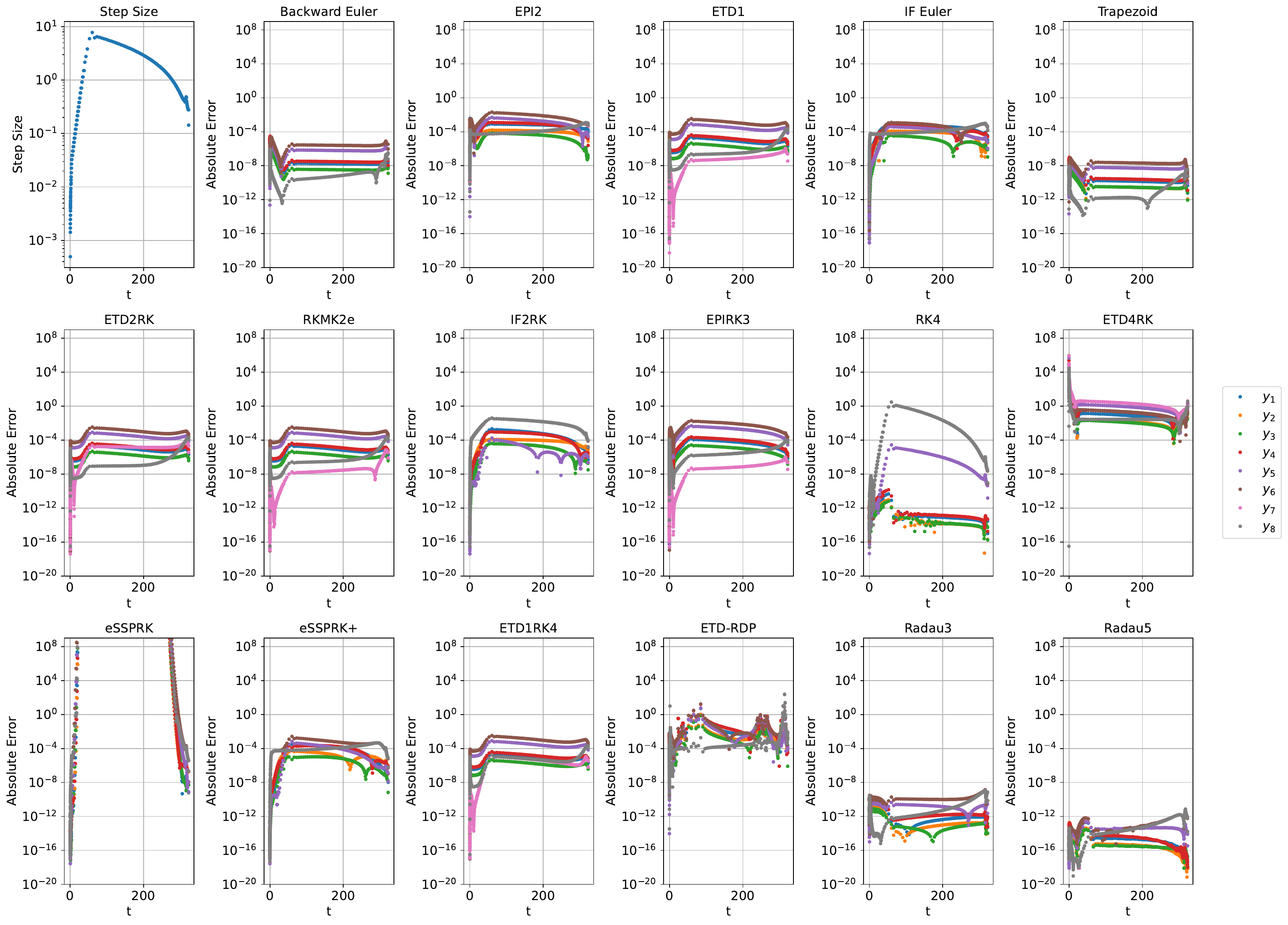}
    \caption{Absolute error for various integration schemes on the HIRES model. The step size is also shown for reference. The data shown correspond to $n=857$ data points across the integration time interval.}
    \label{fig:857_HIRES_plot}
\end{figure*}
\end{turnpage}

\begin{turnpage}
    \begin{figure*}
    \centering
    \includegraphics[width=0.95\linewidth]{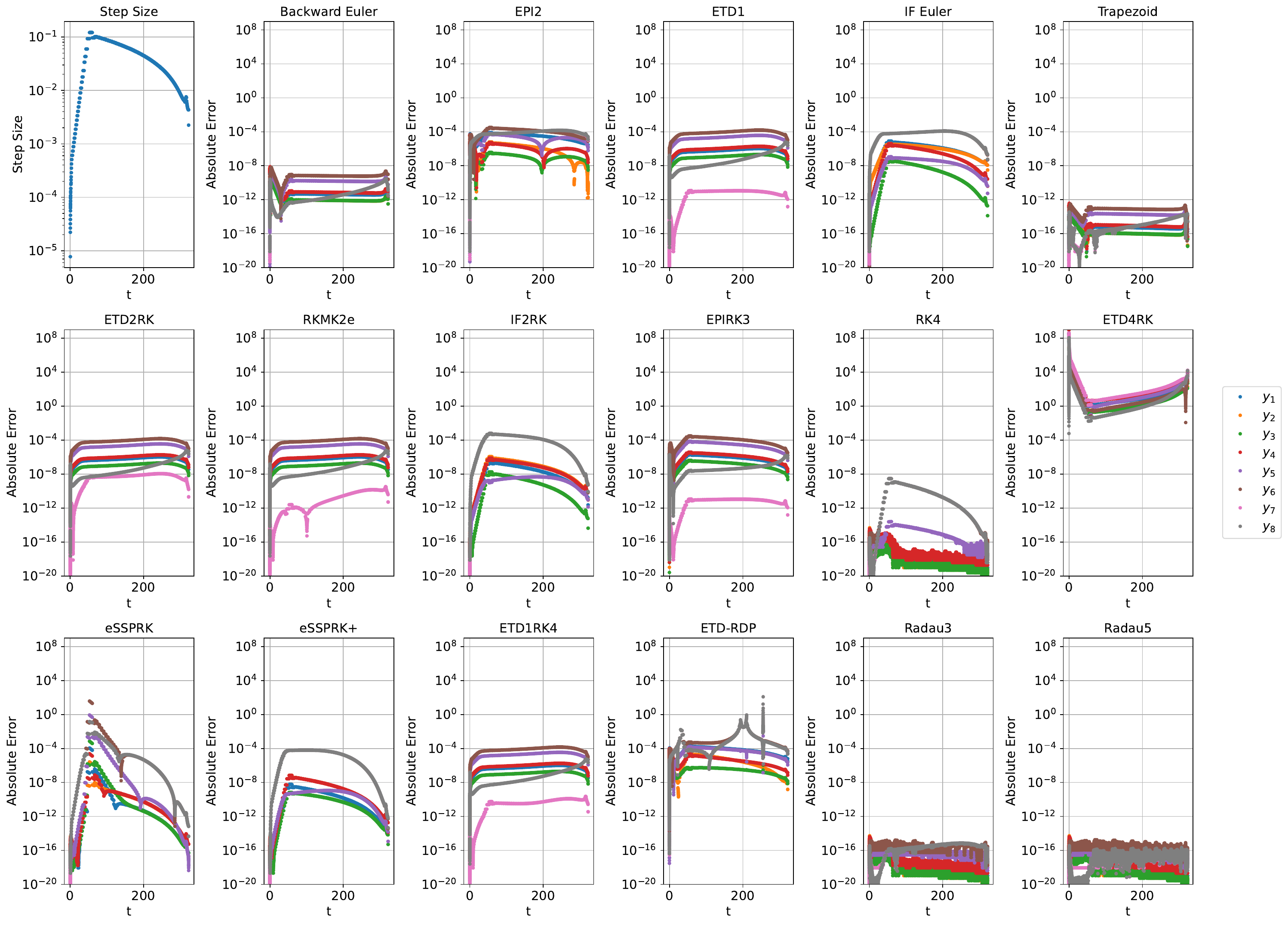}
    \caption{Absolute error for various integration schemes on the HIRES model. The step size is also shown for reference. The data shown correspond to $n=54785$ data points across the integration time interval.}
    \label{fig:54785_HIRES_plot}
\end{figure*}
\end{turnpage}

Figures \ref{fig:56_HIRES_plot}, \ref{fig:857_HIRES_plot}, and \ref{fig:54785_HIRES_plot} show the absolute error for the various integration schemes corresponding to 56, 857, and 54875 data points, respectively. For larger time steps ($n=56$), the eSSPRK, RK4, and ETD4RK methods are the most unstable methods.  eSSPRK is significantly less stable than RK4 and ETD4RK is slightly more stable than RK4. eSSPRK+ is a stable method, whereas eSSPRK is unstable. ETD-RDP, eSSPRK+, EPIRK3, IF2RK, and EPI2 are the methods that perform worse than backward Euler but are still stable methods. Surprisingly, the first-order IF Euler method achieves the highest accuracy, outperforming higher-order methods despite its simplicity and lower computational cost. IF Euler achieves a level of accuracy that lies between the backward Euler and Trapezoid methods. When using a larger step size, the exponential methods fall far short of the accuracy achieved by the Radau methods.

For moderate time steps ($n=857$), the unstable methods are eSSPRK, RK4, ETD-RDP, and ETD4RK. eSSPRK has the worst stability followed by the RK4 method and the ETD4RK method. The remaining exponential integration methods of ETD1, IF Euler, ETD2RK, RKMK2e, IF2RK, EPIRK3, eSSPRK+, and ETDIRK4 all have about the same accuracy, with a performance slightly worst than backward Euler method.

For small time steps ($n=54875$), RK4 is stable with an accuracy between that of backward Euler and Trapezoid method. However, the ETD4RK, eSSPRK, and ETD-RDP methods still are unreliable with very large error. The remaining exponential integration schemes of EPI2, ETD1, IF Euler, ETD2RK, RKMK2e, IF2RK, EPIRK3, eSSPRK+, and ETDIRK4 all seem to saturate at the same accuracy, which has an absolute error about $10^4$ times greater than both backward Euler and the RK4 method.

\begin{turnpage}
\begin{figure*}
    \centering
    \includegraphics[width=0.9\linewidth]{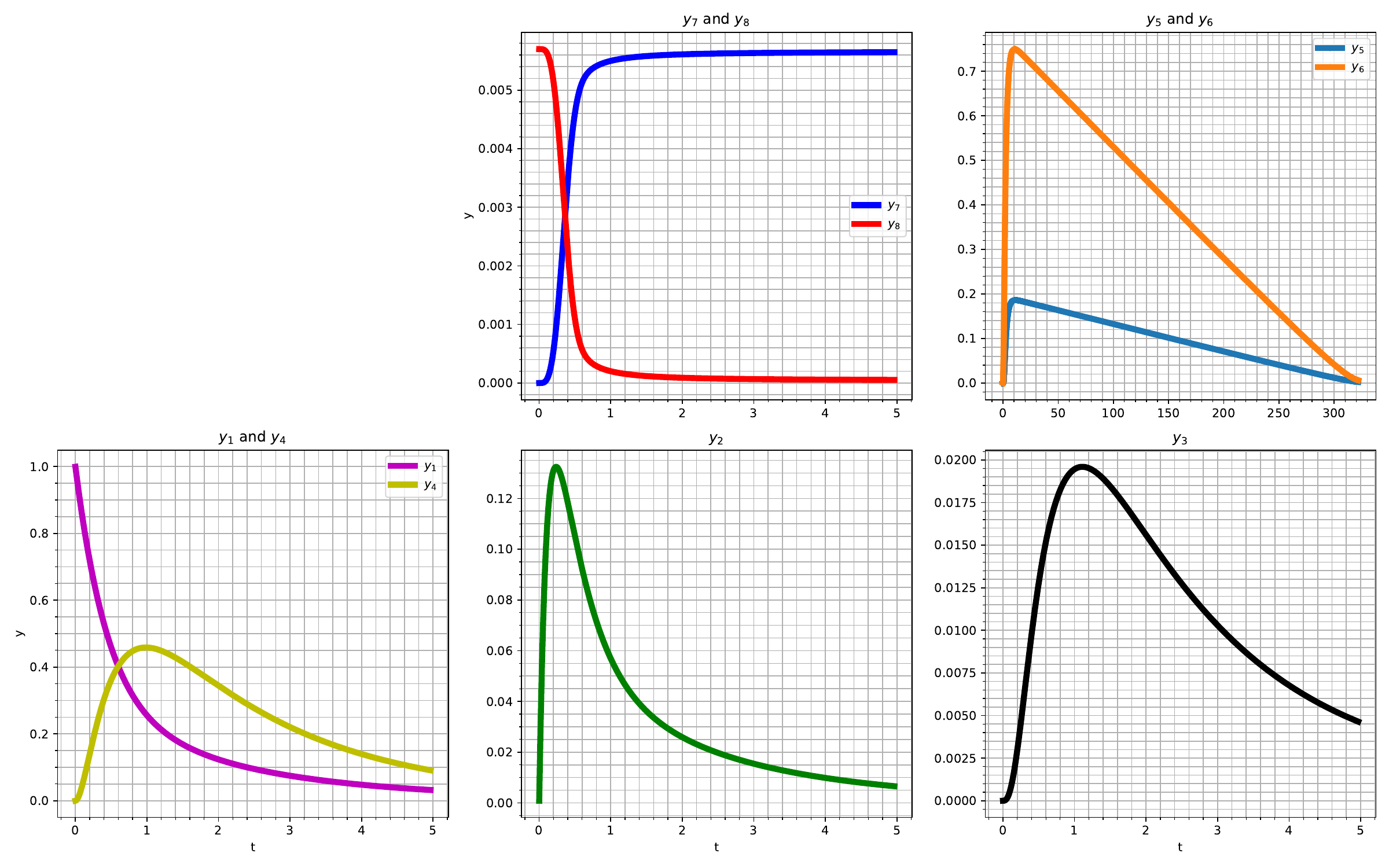}
    \caption{The “High Irradiance RESponse” (HIRES) model}
    \label{fig:HIRES_data}
\end{figure*}
\end{turnpage}

\clearpage
\FloatBarrier
\subsection{Example 3: Robertson Model}

The Robertson model \cite{robertson1966reaction} is given by the following set of equations:

\begin{equation}
\label{eqn:robertson_model}
\begin{aligned}
    &\frac{dy_1}{dt} = -0.04y_1 + 10^4 y_2 y_3, \\
    &\frac{dy_2}{dt} = 0.04y_1 - 10^4 y_2 y_3 - 3 \times 10^7 y_2^2, \\
    &\frac{dy_3}{dt} = 3 \times 10^7 y_2^2, \\
    &y_1(0) = 1, \quad y_2(0)=0, \quad y_3(0)=0, \\
    &t \in [10^{-5}, 10^{7}].
\end{aligned}
\end{equation}

\begin{figure*}
    \centering
    \includegraphics[width=0.9\linewidth]{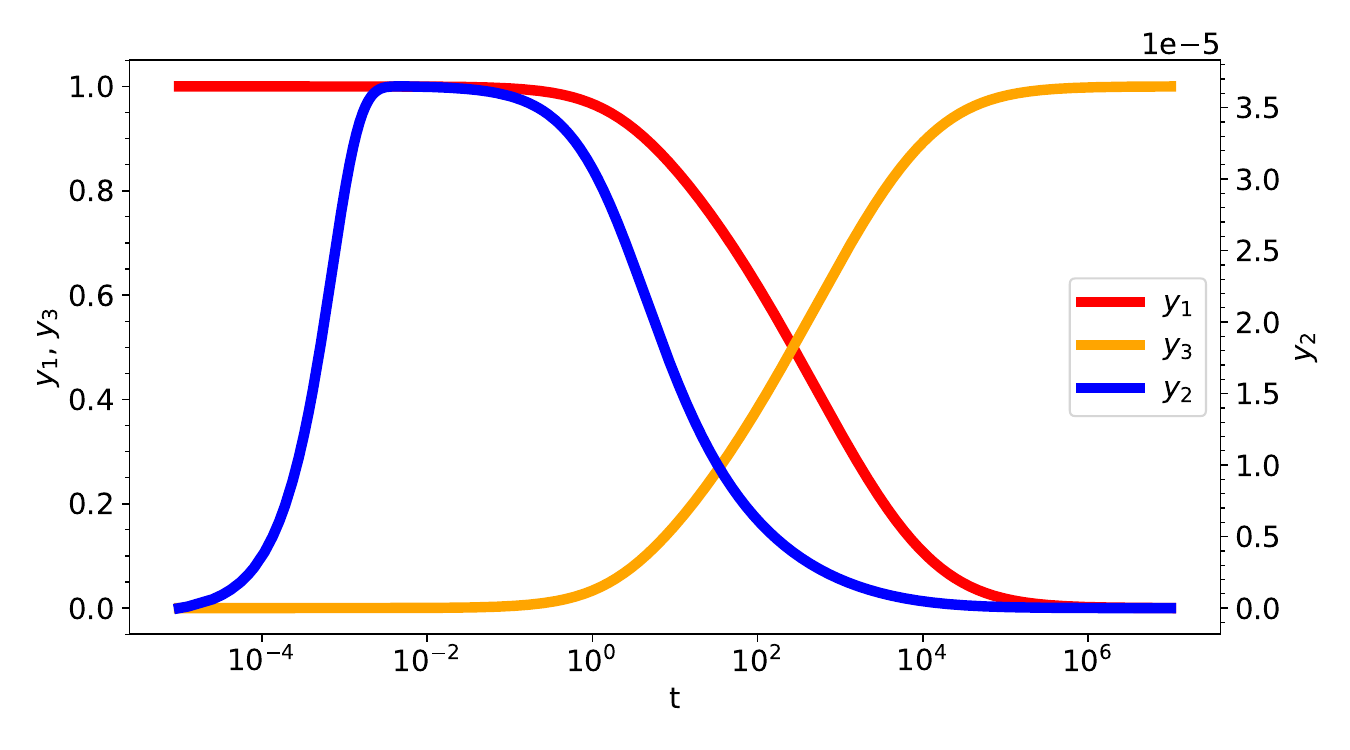}
    \caption{The Robertson model. Due to the different magnitude of the variables, variables $y_1$ and $y_3$ are plotted on the left y-axis and $y_2$ is shown on the right y-axis.}
    \label{fig:robertson_data}
\end{figure*}

\noindent The Robertson model is a set of equations that describe an autocatalytic chemical reaction. Figure \ref{fig:robertson_data} shows the evolution of the chemical species over time. This model has two qualities that make it numerically challenging: the time range is logarithmically spaced over the time interval requiring discretization schemes with stability for really large time steps and one of the species has tiny values around magnitude $10^{-5}$. These difficulties make it an ideal benchmark problem for testing the performance of the discretization schemes.

\begin{turnpage}
    \begin{figure*}
    \centering
    \includegraphics[width=0.95\linewidth]{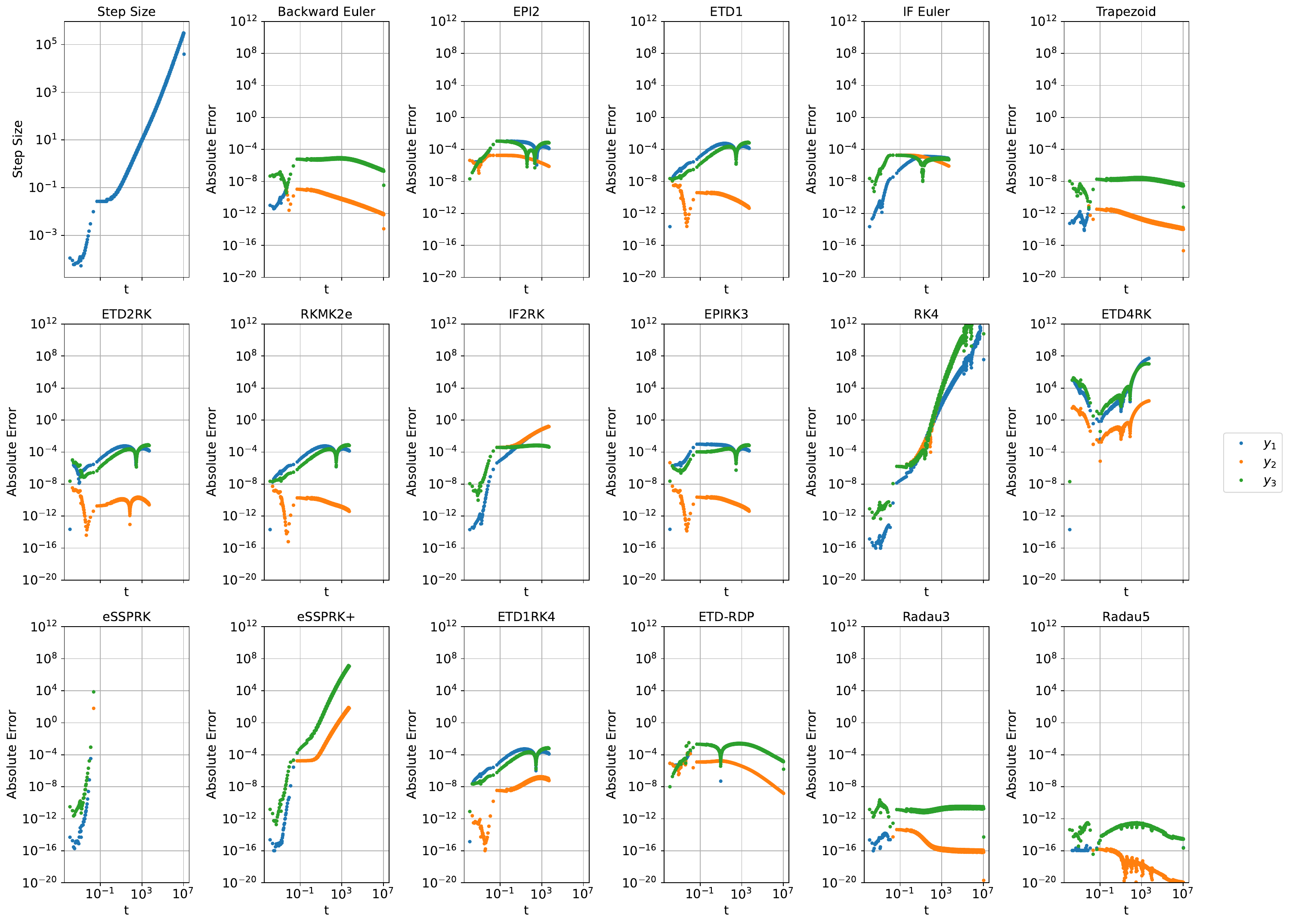}
    \caption{Absolute error for various integration schemes on the Robertson model. The step size is also shown for reference. The data shown correspond to $n=1314$ data points across the integration time interval.}
    \label{fig:1314_Rob_plot}
\end{figure*}
\end{turnpage}

\begin{turnpage}
    \begin{figure*}
    \centering
    \includegraphics[width=0.95\linewidth]{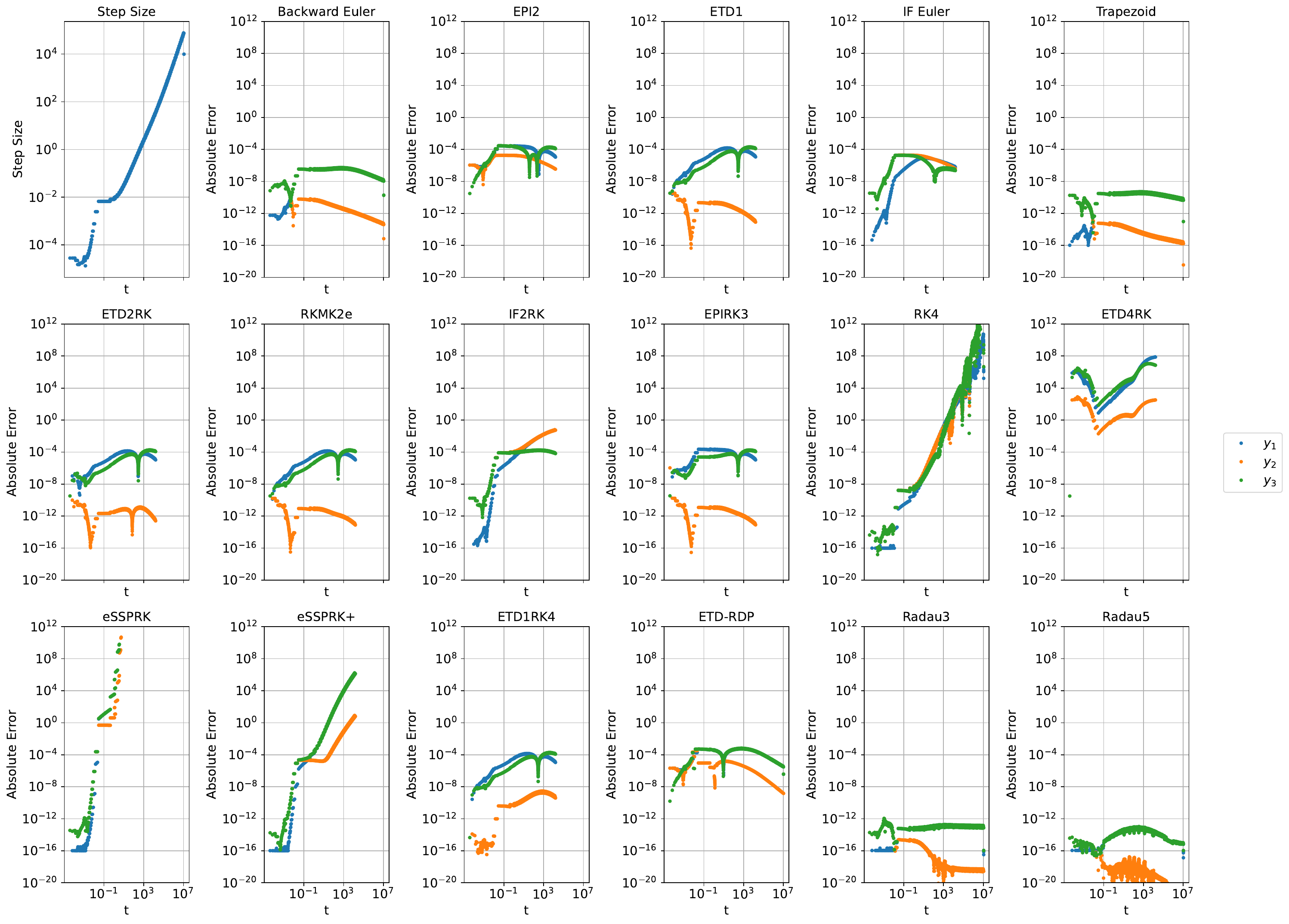}
    \caption{Absolute error for various integration schemes on the Robertson model. The step size is also shown for reference. The data shown correspond to $n=5253$ data points across the integration time interval.}
    \label{fig:5253_Rob_plot}
\end{figure*}
\end{turnpage}

\begin{turnpage}
    \begin{figure*}
    \centering
    \includegraphics[width=0.95\linewidth]{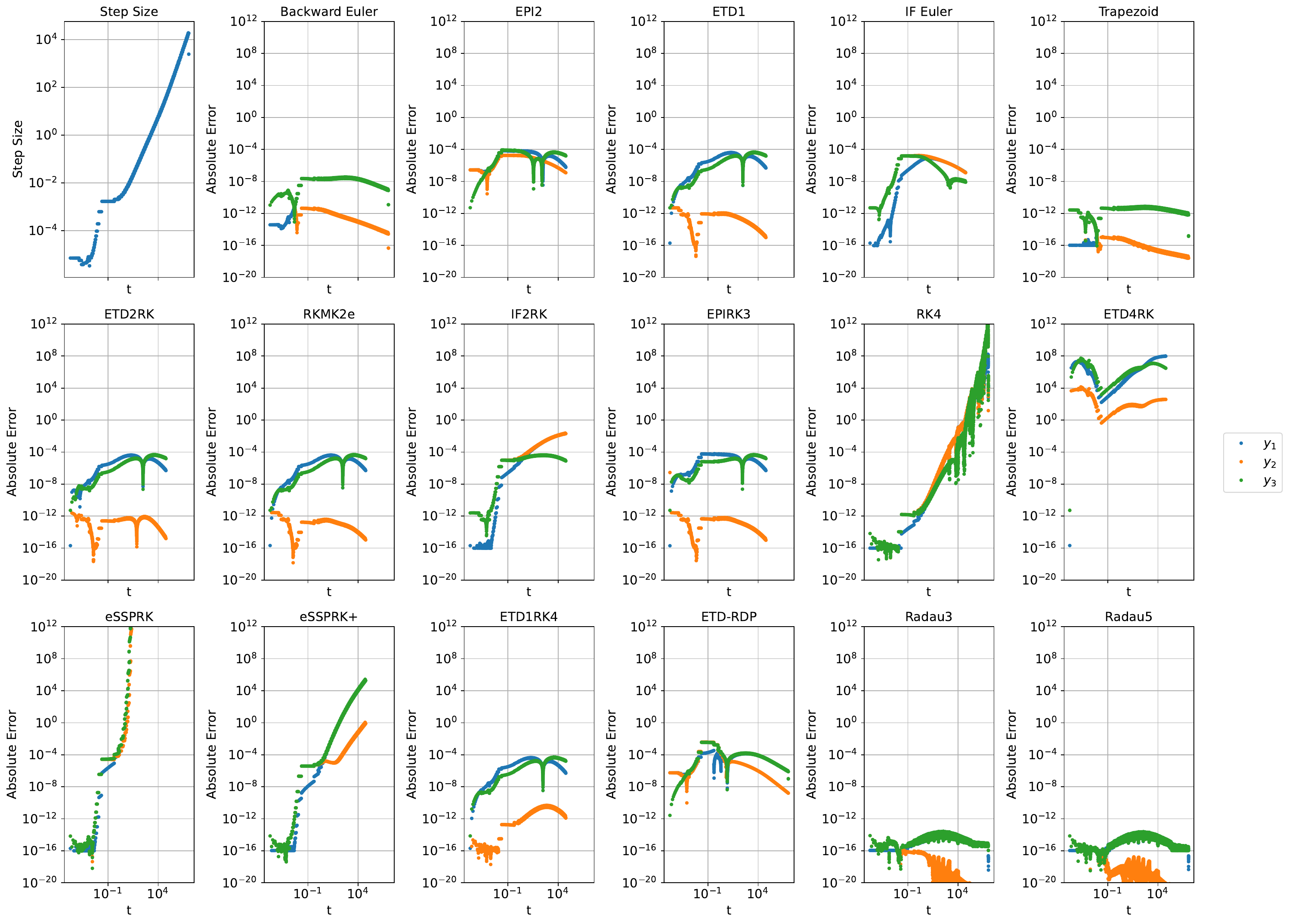}
    \caption{Absolute error for various integration schemes on the Robertson model. The step size is also shown for reference. The data shown correspond to $n=21009$ data points across the integration time interval.}
    \label{fig:21009_Rob_plot}
\end{figure*}
\end{turnpage}

Figures \ref{fig:1314_Rob_plot}, \ref{fig:5253_Rob_plot}, and \ref{fig:21009_Rob_plot} show the absolute error for the various integration schemes corresponding to 1314, 5253, and 21009 data points, respectively. For all of these data points, RK4, ETD4RK, eSSPRK, and eSSPRK+ fail to maintain stability with large integration errors. For extremely large step sizes (larger than $10^3$), the only exponential integration scheme that remains stable is ETD-RDP. The remaining exponential integration methods saturate at an absolute error of around $10^{-4}$, even as the step size is considerably decreased. For reference, for the smallest step size ($n=21009$) backward Euler, Trapezoid method, Radau3, and Radau5 have errors of $10^{-8}$, $10^{-12}$, $10^{-14}$, and $10^{-14}$ respectively. The exponential integration methods lack the accuracy needed to match even the backward Euler method, making them unfavorable for solving the Robertson ODE problem.

\FloatBarrier
\clearpage
\section{Conclusion}

In this study, we evaluated the performance of several promising explicit single-step exponential integration schemes and compared them to classical implicit schemes. Classical implicit integration schemes have been the standard approach for solving stiff ODEs; however, they are slow when applied to applications of ODEs that require repeated integration such as parameter estimation, Bayesian inference, neural ODEs, physics-informed neural networks, and MeshGraphNets. Explicit exponential integration schemes eliminate the need for the nonlinear solver required by implicit schemes, potentially enhancing inference methods and easing the challenges of backpropagating through an ODE solver when implicit schemes are used. Explicit exponential integration methods initially appeared to be a feasible option, but this study has revealed that they do not possess the accuracy required to replace implicit schemes. When the step size was reduced, the backward Euler method was far more accurate than any of the exponential integration schemes. The accuracy of higher-order exponential integration methods does not surpass that of the first-order integrating factor Euler method. Based on this result, we suggest using the IF Euler method for neural ODE applications, as higher-order methods not only lack accuracy but also are computationally more demanding and challenging to backpropagate through. In future work, we confirm the IF Euler’s capability to train stiff neural ODEs, highlighting its viability as a reliable approach for neural differential equations.

\section{Acknowledgements}

Use was made of computational facilities purchased with funds from the National Science Foundation (CNS-1725797) and administered by the Center for Scientific Computing (CSC). The CSC is supported by the California NanoSystems Institute and the Materials Research Science and Engineering Center (MRSEC; NSF DMR 1720256) at UC Santa Barbara. This work was supported in part by NSF awards CNS-1730158, ACI-1540112, ACI-1541349, OAC-1826967, OAC-2112167, CNS-2120019, the University of California Office of the President, and the University of California San Diego's California Institute for Telecommunications and Information Technology/Qualcomm Institute. Thanks to CENIC for the 100Gbps networks. The content of the information does not necessarily reflect the position or the policy of the funding agencies, and no official endorsement should be inferred.  The funders had no role in study design, data collection and analysis, decision to publish, or preparation of the manuscript.

\FloatBarrier
\bibliography{main}

\begin{thebibliography}{75}%
\makeatletter
\providecommand \@ifxundefined [1]{%
 \@ifx{#1\undefined}
}%
\providecommand \@ifnum [1]{%
 \ifnum #1\expandafter \@firstoftwo
 \else \expandafter \@secondoftwo
 \fi
}%
\providecommand \@ifx [1]{%
 \ifx #1\expandafter \@firstoftwo
 \else \expandafter \@secondoftwo
 \fi
}%
\providecommand \natexlab [1]{#1}%
\providecommand \enquote  [1]{``#1''}%
\providecommand \bibnamefont  [1]{#1}%
\providecommand \bibfnamefont [1]{#1}%
\providecommand \citenamefont [1]{#1}%
\providecommand \href@noop [0]{\@secondoftwo}%
\providecommand \href [0]{\begingroup \@sanitize@url \@href}%
\providecommand \@href[1]{\@@startlink{#1}\@@href}%
\providecommand \@@href[1]{\endgroup#1\@@endlink}%
\providecommand \@sanitize@url [0]{\catcode `\\12\catcode `\$12\catcode `\&12\catcode `\#12\catcode `\^12\catcode `\_12\catcode `\%12\relax}%
\providecommand \@@startlink[1]{}%
\providecommand \@@endlink[0]{}%
\providecommand \url  [0]{\begingroup\@sanitize@url \@url }%
\providecommand \@url [1]{\endgroup\@href {#1}{\urlprefix }}%
\providecommand \urlprefix  [0]{URL }%
\providecommand \Eprint [0]{\href }%
\providecommand \doibase [0]{http://dx.doi.org/}%
\providecommand \selectlanguage [0]{\@gobble}%
\providecommand \bibinfo  [0]{\@secondoftwo}%
\providecommand \bibfield  [0]{\@secondoftwo}%
\providecommand \translation [1]{[#1]}%
\providecommand \BibitemOpen [0]{}%
\providecommand \bibitemStop [0]{}%
\providecommand \bibitemNoStop [0]{.\EOS\space}%
\providecommand \EOS [0]{\spacefactor3000\relax}%
\providecommand \BibitemShut  [1]{\csname bibitem#1\endcsname}%
\let\auto@bib@innerbib\@empty
\bibitem [{\citenamefont {Ascher}\ and\ \citenamefont {Petzold}(1998)}]{ascher1998computer}%
  \BibitemOpen
  \bibfield  {author} {\bibinfo {author} {\bibfnamefont {U.~M.}\ \bibnamefont {Ascher}}\ and\ \bibinfo {author} {\bibfnamefont {L.~R.}\ \bibnamefont {Petzold}},\ }\href@noop {} {\emph {\bibinfo {title} {Computer methods for ordinary differential equations and differential-algebraic equations}}},\ Vol.~\bibinfo {volume} {61}\ (\bibinfo  {publisher} {Siam},\ \bibinfo {year} {1998})\BibitemShut {NoStop}%
\bibitem [{\citenamefont {Chen}, \citenamefont {Shao},\ and\ \citenamefont {Ibrahim}(2012)}]{chen2012monte}%
  \BibitemOpen
  \bibfield  {author} {\bibinfo {author} {\bibfnamefont {M.}~\bibnamefont {Chen}}, \bibinfo {author} {\bibfnamefont {Q.}~\bibnamefont {Shao}}, \ and\ \bibinfo {author} {\bibfnamefont {J.}~\bibnamefont {Ibrahim}},\ }\href {https://books.google.com/books?id=4IrbBwAAQBAJ} {\emph {\bibinfo {title} {Monte Carlo Methods in Bayesian Computation}}},\ Springer Series in Statistics\ (\bibinfo  {publisher} {Springer New York},\ \bibinfo {year} {2012})\BibitemShut {NoStop}%
\bibitem [{\citenamefont {Liang}, \citenamefont {Liu},\ and\ \citenamefont {Carroll}(2011)}]{liang2011advanced}%
  \BibitemOpen
  \bibfield  {author} {\bibinfo {author} {\bibfnamefont {F.}~\bibnamefont {Liang}}, \bibinfo {author} {\bibfnamefont {C.}~\bibnamefont {Liu}}, \ and\ \bibinfo {author} {\bibfnamefont {R.}~\bibnamefont {Carroll}},\ }\href {https://books.google.com/books?id=ZmKgUO2PVpIC} {\emph {\bibinfo {title} {Advanced Markov Chain Monte Carlo Methods: Learning from Past Samples}}},\ Wiley Series in Computational Statistics\ (\bibinfo  {publisher} {Wiley},\ \bibinfo {year} {2011})\BibitemShut {NoStop}%
\bibitem [{\citenamefont {McElreath}(2018)}]{mcelreath2018statistical}%
  \BibitemOpen
  \bibfield  {author} {\bibinfo {author} {\bibfnamefont {R.}~\bibnamefont {McElreath}},\ }\href {https://books.google.com/books?id=T3FQDwAAQBAJ} {\emph {\bibinfo {title} {Statistical Rethinking: A Bayesian Course with Examples in R and Stan}}},\ Chapman \& Hall/CRC Texts in Statistical Science\ (\bibinfo  {publisher} {CRC Press},\ \bibinfo {year} {2018})\BibitemShut {NoStop}%
\bibitem [{\citenamefont {Chen}\ \emph {et~al.}(2018)\citenamefont {Chen}, \citenamefont {Rubanova}, \citenamefont {Bettencourt},\ and\ \citenamefont {Duvenaud}}]{chen2018neural}%
  \BibitemOpen
  \bibfield  {author} {\bibinfo {author} {\bibfnamefont {R.~T.}\ \bibnamefont {Chen}}, \bibinfo {author} {\bibfnamefont {Y.}~\bibnamefont {Rubanova}}, \bibinfo {author} {\bibfnamefont {J.}~\bibnamefont {Bettencourt}}, \ and\ \bibinfo {author} {\bibfnamefont {D.~K.}\ \bibnamefont {Duvenaud}},\ }\bibfield  {title} {\enquote {\bibinfo {title} {Neural ordinary differential equations},}\ }\href@noop {} {\bibfield  {journal} {\bibinfo  {journal} {Advances in neural information processing systems}\ }\textbf {\bibinfo {volume} {31}} (\bibinfo {year} {2018})}\BibitemShut {NoStop}%
\bibitem [{\citenamefont {Rubanova}, \citenamefont {Chen},\ and\ \citenamefont {Duvenaud}(2019)}]{latent_ODEs}%
  \BibitemOpen
  \bibfield  {author} {\bibinfo {author} {\bibfnamefont {Y.}~\bibnamefont {Rubanova}}, \bibinfo {author} {\bibfnamefont {R.~T.~Q.}\ \bibnamefont {Chen}}, \ and\ \bibinfo {author} {\bibfnamefont {D.~K.}\ \bibnamefont {Duvenaud}},\ }\bibfield  {title} {\enquote {\bibinfo {title} {Latent ordinary differential equations for irregularly-sampled time series},}\ }in\ \href@noop {} {\emph {\bibinfo {booktitle} {Advances in Neural Information Processing Systems}}},\ Vol.~\bibinfo {volume} {32},\ \bibinfo {editor} {edited by\ \bibinfo {editor} {\bibfnamefont {H.}~\bibnamefont {Wallach}}, \bibinfo {editor} {\bibfnamefont {H.}~\bibnamefont {Larochelle}}, \bibinfo {editor} {\bibfnamefont {A.}~\bibnamefont {Beygelzimer}}, \bibinfo {editor} {\bibfnamefont {F.}~\bibnamefont {d\textquotesingle Alch\'{e}-Buc}}, \bibinfo {editor} {\bibfnamefont {E.}~\bibnamefont {Fox}}, \ and\ \bibinfo {editor} {\bibfnamefont {R.}~\bibnamefont {Garnett}}}\ (\bibinfo  {publisher} {Curran Associates, Inc.},\ \bibinfo {year} {2019})\BibitemShut
  {NoStop}%
\bibitem [{\citenamefont {Dandekar}\ \emph {et~al.}(2020)\citenamefont {Dandekar}, \citenamefont {Dixit}, \citenamefont {Tarek}, \citenamefont {Garcia{-}Valadez},\ and\ \citenamefont {Rackauckas}}]{bayesianneuralode}%
  \BibitemOpen
  \bibfield  {author} {\bibinfo {author} {\bibfnamefont {R.}~\bibnamefont {Dandekar}}, \bibinfo {author} {\bibfnamefont {V.}~\bibnamefont {Dixit}}, \bibinfo {author} {\bibfnamefont {M.}~\bibnamefont {Tarek}}, \bibinfo {author} {\bibfnamefont {A.}~\bibnamefont {Garcia{-}Valadez}}, \ and\ \bibinfo {author} {\bibfnamefont {C.}~\bibnamefont {Rackauckas}},\ }\bibfield  {title} {\enquote {\bibinfo {title} {Bayesian neural ordinary differential equations},}\ }\href@noop {} {\bibfield  {journal} {\bibinfo  {journal} {CoRR}\ }\textbf {\bibinfo {volume} {abs/2012.07244}} (\bibinfo {year} {2020})},\ \Eprint {http://arxiv.org/abs/2012.07244} {2012.07244} \BibitemShut {NoStop}%
\bibitem [{\citenamefont {Li}\ \emph {et~al.}(2020)\citenamefont {Li}, \citenamefont {Wong}, \citenamefont {Chen},\ and\ \citenamefont {Duvenaud}}]{stochastic_neural_ode}%
  \BibitemOpen
  \bibfield  {author} {\bibinfo {author} {\bibfnamefont {X.}~\bibnamefont {Li}}, \bibinfo {author} {\bibfnamefont {T.-K.~L.}\ \bibnamefont {Wong}}, \bibinfo {author} {\bibfnamefont {R.~T.~Q.}\ \bibnamefont {Chen}}, \ and\ \bibinfo {author} {\bibfnamefont {D.}~\bibnamefont {Duvenaud}},\ }\bibfield  {title} {\enquote {\bibinfo {title} {Scalable gradients for stochastic differential equations},}\ }in\ \href@noop {} {\emph {\bibinfo {booktitle} {Proceedings of the Twenty Third International Conference on Artificial Intelligence and Statistics}}},\ \bibinfo {series} {Proceedings of Machine Learning Research}, Vol.\ \bibinfo {volume} {108},\ \bibinfo {editor} {edited by\ \bibinfo {editor} {\bibfnamefont {S.}~\bibnamefont {Chiappa}}\ and\ \bibinfo {editor} {\bibfnamefont {R.}~\bibnamefont {Calandra}}}\ (\bibinfo  {publisher} {PMLR},\ \bibinfo {year} {2020})\ pp.\ \bibinfo {pages} {3870--3882}\BibitemShut {NoStop}%
\bibitem [{\citenamefont {Kidger}\ \emph {et~al.}(2020)\citenamefont {Kidger}, \citenamefont {Morrill}, \citenamefont {Foster},\ and\ \citenamefont {Lyons}}]{kidger2020neural}%
  \BibitemOpen
  \bibfield  {author} {\bibinfo {author} {\bibfnamefont {P.}~\bibnamefont {Kidger}}, \bibinfo {author} {\bibfnamefont {J.}~\bibnamefont {Morrill}}, \bibinfo {author} {\bibfnamefont {J.}~\bibnamefont {Foster}}, \ and\ \bibinfo {author} {\bibfnamefont {T.}~\bibnamefont {Lyons}},\ }\bibfield  {title} {\enquote {\bibinfo {title} {Neural controlled differential equations for irregular time series},}\ }\href@noop {} {\bibfield  {journal} {\bibinfo  {journal} {Advances in Neural Information Processing Systems}\ }\textbf {\bibinfo {volume} {33}},\ \bibinfo {pages} {6696--6707} (\bibinfo {year} {2020})}\BibitemShut {NoStop}%
\bibitem [{\citenamefont {Kidger}(2022)}]{kidger2022neural}%
  \BibitemOpen
  \bibfield  {author} {\bibinfo {author} {\bibfnamefont {P.}~\bibnamefont {Kidger}},\ }\bibfield  {title} {\enquote {\bibinfo {title} {On neural differential equations},}\ }\href@noop {} {\bibfield  {journal} {\bibinfo  {journal} {arXiv preprint arXiv:2202.02435}\ } (\bibinfo {year} {2022})}\BibitemShut {NoStop}%
\bibitem [{\citenamefont {Morrill}\ \emph {et~al.}(2021)\citenamefont {Morrill}, \citenamefont {Salvi}, \citenamefont {Kidger},\ and\ \citenamefont {Foster}}]{morrill2021neural}%
  \BibitemOpen
  \bibfield  {author} {\bibinfo {author} {\bibfnamefont {J.}~\bibnamefont {Morrill}}, \bibinfo {author} {\bibfnamefont {C.}~\bibnamefont {Salvi}}, \bibinfo {author} {\bibfnamefont {P.}~\bibnamefont {Kidger}}, \ and\ \bibinfo {author} {\bibfnamefont {J.}~\bibnamefont {Foster}},\ }\bibfield  {title} {\enquote {\bibinfo {title} {Neural rough differential equations for long time series},}\ }in\ \href@noop {} {\emph {\bibinfo {booktitle} {International Conference on Machine Learning}}}\ (\bibinfo {organization} {PMLR},\ \bibinfo {year} {2021})\ pp.\ \bibinfo {pages} {7829--7838}\BibitemShut {NoStop}%
\bibitem [{\citenamefont {Jia}\ and\ \citenamefont {Benson}(2019)}]{jia2019neural}%
  \BibitemOpen
  \bibfield  {author} {\bibinfo {author} {\bibfnamefont {J.}~\bibnamefont {Jia}}\ and\ \bibinfo {author} {\bibfnamefont {A.~R.}\ \bibnamefont {Benson}},\ }\bibfield  {title} {\enquote {\bibinfo {title} {Neural jump stochastic differential equations},}\ }\href@noop {} {\bibfield  {journal} {\bibinfo  {journal} {Advances in Neural Information Processing Systems}\ }\textbf {\bibinfo {volume} {32}} (\bibinfo {year} {2019})}\BibitemShut {NoStop}%
\bibitem [{\citenamefont {Chen}, \citenamefont {Amos},\ and\ \citenamefont {Nickel}(2020)}]{chen2020learning}%
  \BibitemOpen
  \bibfield  {author} {\bibinfo {author} {\bibfnamefont {R.~T.}\ \bibnamefont {Chen}}, \bibinfo {author} {\bibfnamefont {B.}~\bibnamefont {Amos}}, \ and\ \bibinfo {author} {\bibfnamefont {M.}~\bibnamefont {Nickel}},\ }\bibfield  {title} {\enquote {\bibinfo {title} {Learning neural event functions for ordinary differential equations},}\ }\href@noop {} {\bibfield  {journal} {\bibinfo  {journal} {arXiv preprint arXiv:2011.03902}\ } (\bibinfo {year} {2020})}\BibitemShut {NoStop}%
\bibitem [{\citenamefont {Duvenaud}\ \emph {et~al.}(2023)\citenamefont {Duvenaud}, \citenamefont {Heinonen}, \citenamefont {Tiemann},\ and\ \citenamefont {Welling}}]{dagstuhl}%
  \BibitemOpen
  \bibfield  {author} {\bibinfo {author} {\bibfnamefont {D.}~\bibnamefont {Duvenaud}}, \bibinfo {author} {\bibfnamefont {M.}~\bibnamefont {Heinonen}}, \bibinfo {author} {\bibfnamefont {M.}~\bibnamefont {Tiemann}}, \ and\ \bibinfo {author} {\bibfnamefont {M.}~\bibnamefont {Welling}},\ }\bibfield  {title} {\enquote {\bibinfo {title} {Differential equations and continuous-time deep learning},}\ }\href@noop {} {\bibfield  {journal} {\bibinfo  {journal} {Visualization and Decision Making Design Under Uncertainty}\ ,\ \bibinfo {pages} {19}} (\bibinfo {year} {2023})}\BibitemShut {NoStop}%
\bibitem [{\citenamefont {Fronk}\ and\ \citenamefont {Petzold}(2023)}]{doi:10.1063/5.0130803}%
  \BibitemOpen
  \bibfield  {author} {\bibinfo {author} {\bibfnamefont {C.}~\bibnamefont {Fronk}}\ and\ \bibinfo {author} {\bibfnamefont {L.}~\bibnamefont {Petzold}},\ }\bibfield  {title} {\enquote {\bibinfo {title} {Interpretable polynomial neural ordinary differential equations},}\ }\href {\doibase 10.1063/5.0130803} {\bibfield  {journal} {\bibinfo  {journal} {Chaos: An Interdisciplinary Journal of Nonlinear Science}\ }\textbf {\bibinfo {volume} {33}},\ \bibinfo {pages} {043101} (\bibinfo {year} {2023})}\BibitemShut {NoStop}%
\bibitem [{\citenamefont {Fronk}\ \emph {et~al.}(2024)\citenamefont {Fronk}, \citenamefont {Yun}, \citenamefont {Singh},\ and\ \citenamefont {Petzold}}]{fronk2024bayesian}%
  \BibitemOpen
  \bibfield  {author} {\bibinfo {author} {\bibfnamefont {C.}~\bibnamefont {Fronk}}, \bibinfo {author} {\bibfnamefont {J.}~\bibnamefont {Yun}}, \bibinfo {author} {\bibfnamefont {P.}~\bibnamefont {Singh}}, \ and\ \bibinfo {author} {\bibfnamefont {L.}~\bibnamefont {Petzold}},\ }\bibfield  {title} {\enquote {\bibinfo {title} {Bayesian polynomial neural networks and polynomial neural ordinary differential equations},}\ }\href@noop {} {\bibfield  {journal} {\bibinfo  {journal} {PLOS Computational Biology}\ }\textbf {\bibinfo {volume} {20}},\ \bibinfo {pages} {e1012414} (\bibinfo {year} {2024})}\BibitemShut {NoStop}%
\bibitem [{\citenamefont {Fronk}\ and\ \citenamefont {Petzold}(2024)}]{fronk2024trainingstiffneuralordinary}%
  \BibitemOpen
  \bibfield  {author} {\bibinfo {author} {\bibfnamefont {C.}~\bibnamefont {Fronk}}\ and\ \bibinfo {author} {\bibfnamefont {L.}~\bibnamefont {Petzold}},\ }\href {https://arxiv.org/abs/2410.05592} {\enquote {\bibinfo {title} {Training stiff neural ordinary differential equations with implicit single-step methods},}\ } (\bibinfo {year} {2024}),\ \Eprint {http://arxiv.org/abs/2410.05592} {arXiv:2410.05592 [math.NA]} \BibitemShut {NoStop}%
\bibitem [{\citenamefont {Owhadi}(2015)}]{owhadi2015bayesian}%
  \BibitemOpen
  \bibfield  {author} {\bibinfo {author} {\bibfnamefont {H.}~\bibnamefont {Owhadi}},\ }\bibfield  {title} {\enquote {\bibinfo {title} {Bayesian numerical homogenization},}\ }\href@noop {} {\bibfield  {journal} {\bibinfo  {journal} {Multiscale Modeling \& Simulation}\ }\textbf {\bibinfo {volume} {13}},\ \bibinfo {pages} {812--828} (\bibinfo {year} {2015})}\BibitemShut {NoStop}%
\bibitem [{\citenamefont {Raissi}\ and\ \citenamefont {Karniadakis}(2017)}]{hiddenphysics}%
  \BibitemOpen
  \bibfield  {author} {\bibinfo {author} {\bibfnamefont {M.}~\bibnamefont {Raissi}}\ and\ \bibinfo {author} {\bibfnamefont {G.}~\bibnamefont {Karniadakis}},\ }\bibfield  {title} {\enquote {\bibinfo {title} {Hidden physics models: Machine learning of nonlinear partial differential equations},}\ }\href {\doibase 10.1016/j.jcp.2017.11.039} {\bibfield  {journal} {\bibinfo  {journal} {Journal of Computational Physics}\ }\textbf {\bibinfo {volume} {357}} (\bibinfo {year} {2017}),\ 10.1016/j.jcp.2017.11.039}\BibitemShut {NoStop}%
\bibitem [{\citenamefont {Raissi}, \citenamefont {Perdikaris},\ and\ \citenamefont {Karniadakis}(2018)}]{raissi2018numerical}%
  \BibitemOpen
  \bibfield  {author} {\bibinfo {author} {\bibfnamefont {M.}~\bibnamefont {Raissi}}, \bibinfo {author} {\bibfnamefont {P.}~\bibnamefont {Perdikaris}}, \ and\ \bibinfo {author} {\bibfnamefont {G.~E.}\ \bibnamefont {Karniadakis}},\ }\bibfield  {title} {\enquote {\bibinfo {title} {Numerical gaussian processes for time-dependent and nonlinear partial differential equations},}\ }\href@noop {} {\bibfield  {journal} {\bibinfo  {journal} {SIAM Journal on Scientific Computing}\ }\textbf {\bibinfo {volume} {40}},\ \bibinfo {pages} {A172--A198} (\bibinfo {year} {2018})}\BibitemShut {NoStop}%
\bibitem [{\citenamefont {Raissi}, \citenamefont {Perdikaris},\ and\ \citenamefont {Karniadakis}(2017)}]{raissi2017physics}%
  \BibitemOpen
  \bibfield  {author} {\bibinfo {author} {\bibfnamefont {M.}~\bibnamefont {Raissi}}, \bibinfo {author} {\bibfnamefont {P.}~\bibnamefont {Perdikaris}}, \ and\ \bibinfo {author} {\bibfnamefont {G.~E.}\ \bibnamefont {Karniadakis}},\ }\href@noop {} {\enquote {\bibinfo {title} {Physics informed deep learning (part ii): Data-driven discovery of nonlinear partial differential equations},}\ } (\bibinfo {year} {2017}),\ \Eprint {http://arxiv.org/abs/1711.10566} {arXiv:1711.10566 [cs.AI]} \BibitemShut {NoStop}%
\bibitem [{\citenamefont {Raissi}, \citenamefont {Perdikaris},\ and\ \citenamefont {Karniadakis}(2019)}]{osti_1595805}%
  \BibitemOpen
  \bibfield  {author} {\bibinfo {author} {\bibfnamefont {M.}~\bibnamefont {Raissi}}, \bibinfo {author} {\bibfnamefont {P.}~\bibnamefont {Perdikaris}}, \ and\ \bibinfo {author} {\bibfnamefont {G.~E.}\ \bibnamefont {Karniadakis}},\ }\bibfield  {title} {\enquote {\bibinfo {title} {Physics-informed neural networks: A deep learning framework for solving forward and inverse problems involving nonlinear partial differential equations},}\ }\href@noop {} {\bibfield  {journal} {\bibinfo  {journal} {Journal of Computational physics}\ }\textbf {\bibinfo {volume} {378}},\ \bibinfo {pages} {686--707} (\bibinfo {year} {2019})}\BibitemShut {NoStop}%
\bibitem [{\citenamefont {Cuomo}\ \emph {et~al.}(2022)\citenamefont {Cuomo}, \citenamefont {Di~Cola}, \citenamefont {Giampaolo}, \citenamefont {Rozza}, \citenamefont {Raissi},\ and\ \citenamefont {Piccialli}}]{cuomo2022scientific}%
  \BibitemOpen
  \bibfield  {author} {\bibinfo {author} {\bibfnamefont {S.}~\bibnamefont {Cuomo}}, \bibinfo {author} {\bibfnamefont {V.~S.}\ \bibnamefont {Di~Cola}}, \bibinfo {author} {\bibfnamefont {F.}~\bibnamefont {Giampaolo}}, \bibinfo {author} {\bibfnamefont {G.}~\bibnamefont {Rozza}}, \bibinfo {author} {\bibfnamefont {M.}~\bibnamefont {Raissi}}, \ and\ \bibinfo {author} {\bibfnamefont {F.}~\bibnamefont {Piccialli}},\ }\bibfield  {title} {\enquote {\bibinfo {title} {Scientific machine learning through physics--informed neural networks: Where we are and what’s next},}\ }\href@noop {} {\bibfield  {journal} {\bibinfo  {journal} {Journal of Scientific Computing}\ }\textbf {\bibinfo {volume} {92}},\ \bibinfo {pages} {88} (\bibinfo {year} {2022})}\BibitemShut {NoStop}%
\bibitem [{\citenamefont {Cai}\ \emph {et~al.}(2021)\citenamefont {Cai}, \citenamefont {Mao}, \citenamefont {Wang}, \citenamefont {Yin},\ and\ \citenamefont {Karniadakis}}]{cai2021physics}%
  \BibitemOpen
  \bibfield  {author} {\bibinfo {author} {\bibfnamefont {S.}~\bibnamefont {Cai}}, \bibinfo {author} {\bibfnamefont {Z.}~\bibnamefont {Mao}}, \bibinfo {author} {\bibfnamefont {Z.}~\bibnamefont {Wang}}, \bibinfo {author} {\bibfnamefont {M.}~\bibnamefont {Yin}}, \ and\ \bibinfo {author} {\bibfnamefont {G.~E.}\ \bibnamefont {Karniadakis}},\ }\bibfield  {title} {\enquote {\bibinfo {title} {Physics-informed neural networks (pinns) for fluid mechanics: A review},}\ }\href@noop {} {\bibfield  {journal} {\bibinfo  {journal} {Acta Mechanica Sinica}\ }\textbf {\bibinfo {volume} {37}},\ \bibinfo {pages} {1727--1738} (\bibinfo {year} {2021})}\BibitemShut {NoStop}%
\bibitem [{\citenamefont {Pfaff}\ \emph {et~al.}(2021)\citenamefont {Pfaff}, \citenamefont {Fortunato}, \citenamefont {Sanchez-Gonzalez},\ and\ \citenamefont {Battaglia}}]{pfaff2021learningmeshbasedsimulationgraph}%
  \BibitemOpen
  \bibfield  {author} {\bibinfo {author} {\bibfnamefont {T.}~\bibnamefont {Pfaff}}, \bibinfo {author} {\bibfnamefont {M.}~\bibnamefont {Fortunato}}, \bibinfo {author} {\bibfnamefont {A.}~\bibnamefont {Sanchez-Gonzalez}}, \ and\ \bibinfo {author} {\bibfnamefont {P.~W.}\ \bibnamefont {Battaglia}},\ }\href {https://arxiv.org/abs/2010.03409} {\enquote {\bibinfo {title} {Learning mesh-based simulation with graph networks},}\ } (\bibinfo {year} {2021}),\ \Eprint {http://arxiv.org/abs/2010.03409} {arXiv:2010.03409 [cs.LG]} \BibitemShut {NoStop}%
\bibitem [{\citenamefont {Trefethen}(2000)}]{trefethen2000spectral}%
  \BibitemOpen
  \bibfield  {author} {\bibinfo {author} {\bibfnamefont {L.~N.}\ \bibnamefont {Trefethen}},\ }\href@noop {} {\emph {\bibinfo {title} {Spectral methods in MATLAB}}}\ (\bibinfo  {publisher} {SIAM},\ \bibinfo {year} {2000})\BibitemShut {NoStop}%
\bibitem [{\citenamefont {Boyd}(2001)}]{boyd2001chebyshev}%
  \BibitemOpen
  \bibfield  {author} {\bibinfo {author} {\bibfnamefont {J.~P.}\ \bibnamefont {Boyd}},\ }\href@noop {} {\emph {\bibinfo {title} {Chebyshev and Fourier spectral methods}}}\ (\bibinfo  {publisher} {Courier Corporation},\ \bibinfo {year} {2001})\BibitemShut {NoStop}%
\bibitem [{\citenamefont {Canuto}\ \emph {et~al.}(1988)\citenamefont {Canuto}, \citenamefont {Hussaini}, \citenamefont {Quarteroni},\ and\ \citenamefont {Zang}}]{canuto1988spectral}%
  \BibitemOpen
  \bibfield  {author} {\bibinfo {author} {\bibfnamefont {C.}~\bibnamefont {Canuto}}, \bibinfo {author} {\bibfnamefont {M.~Y.}\ \bibnamefont {Hussaini}}, \bibinfo {author} {\bibfnamefont {A.}~\bibnamefont {Quarteroni}}, \ and\ \bibinfo {author} {\bibfnamefont {T.~A.}\ \bibnamefont {Zang}},\ }\href {\doibase 10.1007/978-3-642-84108-8} {\emph {\bibinfo {title} {Spectral Methods in Fluid Dynamics}}},\ \bibinfo {edition} {1st}\ ed.,\ Scientific Computation\ (\bibinfo  {publisher} {Springer-Verlag Berlin Heidelberg},\ \bibinfo {year} {1988})\BibitemShut {NoStop}%
\bibitem [{\citenamefont {Ascher}(2008)}]{ascher2008numerical}%
  \BibitemOpen
  \bibfield  {author} {\bibinfo {author} {\bibfnamefont {U.~M.}\ \bibnamefont {Ascher}},\ }\href@noop {} {\emph {\bibinfo {title} {Numerical methods for evolutionary differential equations}}}\ (\bibinfo  {publisher} {SIAM},\ \bibinfo {year} {2008})\BibitemShut {NoStop}%
\bibitem [{\citenamefont {Miranker}(2001)}]{miranker2001numerical}%
  \BibitemOpen
  \bibfield  {author} {\bibinfo {author} {\bibfnamefont {A.}~\bibnamefont {Miranker}},\ }\href@noop {} {\emph {\bibinfo {title} {Numerical Methods for Stiff Equations and Singular Perturbation Problems: and singular perturbation problems}}},\ Vol.~\bibinfo {volume} {5}\ (\bibinfo  {publisher} {Springer Science \& Business Media},\ \bibinfo {year} {2001})\BibitemShut {NoStop}%
\bibitem [{\citenamefont {Axelsson}(1969)}]{axelsson1969class}%
  \BibitemOpen
  \bibfield  {author} {\bibinfo {author} {\bibfnamefont {O.}~\bibnamefont {Axelsson}},\ }\bibfield  {title} {\enquote {\bibinfo {title} {A class of a-stable methods},}\ }\href@noop {} {\bibfield  {journal} {\bibinfo  {journal} {BIT Numerical Mathematics}\ }\textbf {\bibinfo {volume} {9}},\ \bibinfo {pages} {185--199} (\bibinfo {year} {1969})}\BibitemShut {NoStop}%
\bibitem [{\citenamefont {Ehle}(1969)}]{ehle1969pade}%
  \BibitemOpen
  \bibfield  {author} {\bibinfo {author} {\bibfnamefont {B.~L.}\ \bibnamefont {Ehle}},\ }\emph {\bibinfo {title} {On Pad{\'e} approximations to the exponential function and A-stable methods for the numerical solution of initial value problems}},\ \href@noop {} {Ph.D. thesis},\ \bibinfo  {school} {University of Waterloo Waterloo, Ontario} (\bibinfo {year} {1969})\BibitemShut {NoStop}%
\bibitem [{\citenamefont {Hairer}\ and\ \citenamefont {Wanner}(1999)}]{hairer1999stiff}%
  \BibitemOpen
  \bibfield  {author} {\bibinfo {author} {\bibfnamefont {E.}~\bibnamefont {Hairer}}\ and\ \bibinfo {author} {\bibfnamefont {G.}~\bibnamefont {Wanner}},\ }\bibfield  {title} {\enquote {\bibinfo {title} {Stiff differential equations solved by radau methods},}\ }\href@noop {} {\bibfield  {journal} {\bibinfo  {journal} {Journal of Computational and Applied Mathematics}\ }\textbf {\bibinfo {volume} {111}},\ \bibinfo {pages} {93--111} (\bibinfo {year} {1999})}\BibitemShut {NoStop}%
\bibitem [{\citenamefont {Lawson}(1967)}]{lawson1967generalized}%
  \BibitemOpen
  \bibfield  {author} {\bibinfo {author} {\bibfnamefont {J.~D.}\ \bibnamefont {Lawson}},\ }\bibfield  {title} {\enquote {\bibinfo {title} {Generalized runge-kutta processes for stable systems with large lipschitz constants},}\ }\href@noop {} {\bibfield  {journal} {\bibinfo  {journal} {SIAM Journal on Numerical Analysis}\ }\textbf {\bibinfo {volume} {4}},\ \bibinfo {pages} {372--380} (\bibinfo {year} {1967})}\BibitemShut {NoStop}%
\bibitem [{\citenamefont {Fornberg}\ and\ \citenamefont {Driscoll}(1999)}]{fornberg1999fast}%
  \BibitemOpen
  \bibfield  {author} {\bibinfo {author} {\bibfnamefont {B.}~\bibnamefont {Fornberg}}\ and\ \bibinfo {author} {\bibfnamefont {T.~A.}\ \bibnamefont {Driscoll}},\ }\bibfield  {title} {\enquote {\bibinfo {title} {A fast spectral algorithm for nonlinear wave equations with linear dispersion},}\ }\href@noop {} {\bibfield  {journal} {\bibinfo  {journal} {Journal of Computational Physics}\ }\textbf {\bibinfo {volume} {155}},\ \bibinfo {pages} {456--467} (\bibinfo {year} {1999})}\BibitemShut {NoStop}%
\bibitem [{\citenamefont {Cox}\ and\ \citenamefont {Matthews}(2002)}]{cox2002exponential}%
  \BibitemOpen
  \bibfield  {author} {\bibinfo {author} {\bibfnamefont {S.~M.}\ \bibnamefont {Cox}}\ and\ \bibinfo {author} {\bibfnamefont {P.~C.}\ \bibnamefont {Matthews}},\ }\bibfield  {title} {\enquote {\bibinfo {title} {Exponential time differencing for stiff systems},}\ }\href@noop {} {\bibfield  {journal} {\bibinfo  {journal} {Journal of Computational Physics}\ }\textbf {\bibinfo {volume} {176}},\ \bibinfo {pages} {430--455} (\bibinfo {year} {2002})}\BibitemShut {NoStop}%
\bibitem [{\citenamefont {Beylkin}, \citenamefont {Keiser},\ and\ \citenamefont {Vozovoi}(1998)}]{beylkin1998new}%
  \BibitemOpen
  \bibfield  {author} {\bibinfo {author} {\bibfnamefont {G.}~\bibnamefont {Beylkin}}, \bibinfo {author} {\bibfnamefont {J.~M.}\ \bibnamefont {Keiser}}, \ and\ \bibinfo {author} {\bibfnamefont {L.}~\bibnamefont {Vozovoi}},\ }\bibfield  {title} {\enquote {\bibinfo {title} {A new class of time discretization schemes for the solution of nonlinear pdes},}\ }\href@noop {} {\bibfield  {journal} {\bibinfo  {journal} {Journal of computational physics}\ }\textbf {\bibinfo {volume} {147}},\ \bibinfo {pages} {362--387} (\bibinfo {year} {1998})}\BibitemShut {NoStop}%
\bibitem [{\citenamefont {Certaine}(1960)}]{certaine1960solution}%
  \BibitemOpen
  \bibfield  {author} {\bibinfo {author} {\bibfnamefont {J.}~\bibnamefont {Certaine}},\ }\bibfield  {title} {\enquote {\bibinfo {title} {The solution of ordinary differential equations with large time constants},}\ }in\ \href@noop {} {\emph {\bibinfo {booktitle} {Mathematical Methods for Digital Computers}}}\ (\bibinfo  {publisher} {Wiley},\ \bibinfo {address} {New York},\ \bibinfo {year} {1960})\ pp.\ \bibinfo {pages} {128--132}\BibitemShut {NoStop}%
\bibitem [{\citenamefont {Hochbruck}\ and\ \citenamefont {Ostermann}(2010)}]{hochbruck2010exponential}%
  \BibitemOpen
  \bibfield  {author} {\bibinfo {author} {\bibfnamefont {M.}~\bibnamefont {Hochbruck}}\ and\ \bibinfo {author} {\bibfnamefont {A.}~\bibnamefont {Ostermann}},\ }\bibfield  {title} {\enquote {\bibinfo {title} {Exponential integrators},}\ }\href@noop {} {\bibfield  {journal} {\bibinfo  {journal} {Acta Numerica}\ }\textbf {\bibinfo {volume} {19}},\ \bibinfo {pages} {209--286} (\bibinfo {year} {2010})}\BibitemShut {NoStop}%
\bibitem [{\citenamefont {Verwer}\ and\ \citenamefont {Van~Loon}(1994)}]{verwer1994evaluation}%
  \BibitemOpen
  \bibfield  {author} {\bibinfo {author} {\bibfnamefont {J.}~\bibnamefont {Verwer}}\ and\ \bibinfo {author} {\bibfnamefont {M.}~\bibnamefont {Van~Loon}},\ }\bibfield  {title} {\enquote {\bibinfo {title} {An evaluation of explicit pseudo-steady-state approximation schemes for stiff ode systems from chemical kinetics},}\ }\href@noop {} {\bibfield  {journal} {\bibinfo  {journal} {Journal of computational physics (Print)}\ }\textbf {\bibinfo {volume} {113}},\ \bibinfo {pages} {347--352} (\bibinfo {year} {1994})}\BibitemShut {NoStop}%
\bibitem [{\citenamefont {Maset}\ and\ \citenamefont {Zennaro}(2009)}]{maset2009unconditional}%
  \BibitemOpen
  \bibfield  {author} {\bibinfo {author} {\bibfnamefont {S.}~\bibnamefont {Maset}}\ and\ \bibinfo {author} {\bibfnamefont {M.}~\bibnamefont {Zennaro}},\ }\bibfield  {title} {\enquote {\bibinfo {title} {Unconditional stability of explicit exponential runge-kutta methods for semi-linear ordinary differential equations},}\ }\href@noop {} {\bibfield  {journal} {\bibinfo  {journal} {Mathematics of computation}\ }\textbf {\bibinfo {volume} {78}},\ \bibinfo {pages} {957--967} (\bibinfo {year} {2009})}\BibitemShut {NoStop}%
\bibitem [{\citenamefont {Iyiola}\ and\ \citenamefont {Wade}(2018)}]{iyiola2018exponential}%
  \BibitemOpen
  \bibfield  {author} {\bibinfo {author} {\bibfnamefont {O.~S.}\ \bibnamefont {Iyiola}}\ and\ \bibinfo {author} {\bibfnamefont {B.~A.}\ \bibnamefont {Wade}},\ }\bibfield  {title} {\enquote {\bibinfo {title} {Exponential integrator methods for systems of non-linear space-fractional models with super-diffusion processes in pattern formation},}\ }\href@noop {} {\bibfield  {journal} {\bibinfo  {journal} {Computers \& Mathematics with Applications}\ }\textbf {\bibinfo {volume} {75}},\ \bibinfo {pages} {3719--3736} (\bibinfo {year} {2018})}\BibitemShut {NoStop}%
\bibitem [{\citenamefont {Tokman}(2006)}]{tokman2006efficient}%
  \BibitemOpen
  \bibfield  {author} {\bibinfo {author} {\bibfnamefont {M.}~\bibnamefont {Tokman}},\ }\bibfield  {title} {\enquote {\bibinfo {title} {Efficient integration of large stiff systems of odes with exponential propagation iterative (epi) methods},}\ }\href@noop {} {\bibfield  {journal} {\bibinfo  {journal} {Journal of Computational Physics}\ }\textbf {\bibinfo {volume} {213}},\ \bibinfo {pages} {748--776} (\bibinfo {year} {2006})}\BibitemShut {NoStop}%
\bibitem [{\citenamefont {Krogstad}(2005)}]{krogstad2005generalized}%
  \BibitemOpen
  \bibfield  {author} {\bibinfo {author} {\bibfnamefont {S.}~\bibnamefont {Krogstad}},\ }\bibfield  {title} {\enquote {\bibinfo {title} {Generalized integrating factor methods for stiff pdes},}\ }\href@noop {} {\bibfield  {journal} {\bibinfo  {journal} {Journal of Computational Physics}\ }\textbf {\bibinfo {volume} {203}},\ \bibinfo {pages} {72--88} (\bibinfo {year} {2005})}\BibitemShut {NoStop}%
\bibitem [{\citenamefont {Li}\ and\ \citenamefont {Li}(2023)}]{li2023low}%
  \BibitemOpen
  \bibfield  {author} {\bibinfo {author} {\bibfnamefont {J.}~\bibnamefont {Li}}\ and\ \bibinfo {author} {\bibfnamefont {D.}~\bibnamefont {Li}},\ }\bibfield  {title} {\enquote {\bibinfo {title} {Low-rank exponential integrators for differential riccati equation},}\ }in\ \href@noop {} {\emph {\bibinfo {booktitle} {Fuzzy Systems and Data Mining IX}}}\ (\bibinfo  {publisher} {IOS Press},\ \bibinfo {year} {2023})\ pp.\ \bibinfo {pages} {372--388}\BibitemShut {NoStop}%
\bibitem [{\citenamefont {Maday}, \citenamefont {Patera},\ and\ \citenamefont {R{\o}nquist}(1990)}]{maday1990operator}%
  \BibitemOpen
  \bibfield  {author} {\bibinfo {author} {\bibfnamefont {Y.}~\bibnamefont {Maday}}, \bibinfo {author} {\bibfnamefont {A.~T.}\ \bibnamefont {Patera}}, \ and\ \bibinfo {author} {\bibfnamefont {E.~M.}\ \bibnamefont {R{\o}nquist}},\ }\bibfield  {title} {\enquote {\bibinfo {title} {An operator-integration-factor splitting method for time-dependent problems: application to incompressible fluid flow},}\ }\href@noop {} {\bibfield  {journal} {\bibinfo  {journal} {Journal of Scientific Computing}\ }\textbf {\bibinfo {volume} {5}},\ \bibinfo {pages} {263--292} (\bibinfo {year} {1990})}\BibitemShut {NoStop}%
\bibitem [{\citenamefont {Isherwood}, \citenamefont {Grant},\ and\ \citenamefont {Gottlieb}(2019)}]{isherwood2019strong}%
  \BibitemOpen
  \bibfield  {author} {\bibinfo {author} {\bibfnamefont {L.}~\bibnamefont {Isherwood}}, \bibinfo {author} {\bibfnamefont {Z.~J.}\ \bibnamefont {Grant}}, \ and\ \bibinfo {author} {\bibfnamefont {S.}~\bibnamefont {Gottlieb}},\ }\bibfield  {title} {\enquote {\bibinfo {title} {Strong stability preserving integrating factor two-step runge--kutta methods},}\ }\href@noop {} {\bibfield  {journal} {\bibinfo  {journal} {Journal of Scientific Computing}\ }\textbf {\bibinfo {volume} {81}},\ \bibinfo {pages} {1446--1471} (\bibinfo {year} {2019})}\BibitemShut {NoStop}%
\bibitem [{\citenamefont {Carrillo}\ \emph {et~al.}(2003)\citenamefont {Carrillo}, \citenamefont {Gamba}, \citenamefont {Majorana},\ and\ \citenamefont {Shu}}]{carrillo2003weno}%
  \BibitemOpen
  \bibfield  {author} {\bibinfo {author} {\bibfnamefont {J.~A.}\ \bibnamefont {Carrillo}}, \bibinfo {author} {\bibfnamefont {I.~M.}\ \bibnamefont {Gamba}}, \bibinfo {author} {\bibfnamefont {A.}~\bibnamefont {Majorana}}, \ and\ \bibinfo {author} {\bibfnamefont {C.-W.}\ \bibnamefont {Shu}},\ }\bibfield  {title} {\enquote {\bibinfo {title} {A weno-solver for the transients of boltzmann--poisson system for semiconductor devices: performance and comparisons with monte carlo methods},}\ }\href@noop {} {\bibfield  {journal} {\bibinfo  {journal} {Journal of Computational Physics}\ }\textbf {\bibinfo {volume} {184}},\ \bibinfo {pages} {498--525} (\bibinfo {year} {2003})}\BibitemShut {NoStop}%
\bibitem [{\citenamefont {Shu}(1988)}]{shu1988total}%
  \BibitemOpen
  \bibfield  {author} {\bibinfo {author} {\bibfnamefont {C.-W.}\ \bibnamefont {Shu}},\ }\bibfield  {title} {\enquote {\bibinfo {title} {Total-variation-diminishing time discretizations},}\ }\href@noop {} {\bibfield  {journal} {\bibinfo  {journal} {SIAM Journal on Scientific and Statistical Computing}\ }\textbf {\bibinfo {volume} {9}},\ \bibinfo {pages} {1073--1084} (\bibinfo {year} {1988})}\BibitemShut {NoStop}%
\bibitem [{\citenamefont {Shu}\ and\ \citenamefont {Osher}(1988)}]{shu1988efficient}%
  \BibitemOpen
  \bibfield  {author} {\bibinfo {author} {\bibfnamefont {C.-W.}\ \bibnamefont {Shu}}\ and\ \bibinfo {author} {\bibfnamefont {S.}~\bibnamefont {Osher}},\ }\bibfield  {title} {\enquote {\bibinfo {title} {Efficient implementation of essentially non-oscillatory shock-capturing schemes},}\ }\href@noop {} {\bibfield  {journal} {\bibinfo  {journal} {Journal of computational physics}\ }\textbf {\bibinfo {volume} {77}},\ \bibinfo {pages} {439--471} (\bibinfo {year} {1988})}\BibitemShut {NoStop}%
\bibitem [{\citenamefont {Ferracina}\ and\ \citenamefont {Spijker}(2004)}]{ferracina2004stepsize}%
  \BibitemOpen
  \bibfield  {author} {\bibinfo {author} {\bibfnamefont {L.}~\bibnamefont {Ferracina}}\ and\ \bibinfo {author} {\bibfnamefont {M.~N.}\ \bibnamefont {Spijker}},\ }\bibfield  {title} {\enquote {\bibinfo {title} {Stepsize restrictions for the total-variation-diminishing property in general runge--kutta methods},}\ }\href@noop {} {\bibfield  {journal} {\bibinfo  {journal} {SIAM journal on numerical analysis}\ }\textbf {\bibinfo {volume} {42}},\ \bibinfo {pages} {1073--1093} (\bibinfo {year} {2004})}\BibitemShut {NoStop}%
\bibitem [{\citenamefont {Ferracina}\ and\ \citenamefont {Spijker}(2005)}]{ferracina2005extension}%
  \BibitemOpen
  \bibfield  {author} {\bibinfo {author} {\bibfnamefont {L.}~\bibnamefont {Ferracina}}\ and\ \bibinfo {author} {\bibfnamefont {M.}~\bibnamefont {Spijker}},\ }\bibfield  {title} {\enquote {\bibinfo {title} {An extension and analysis of the shu-osher representation of runge-kutta methods},}\ }\href@noop {} {\bibfield  {journal} {\bibinfo  {journal} {Mathematics of Computation}\ }\textbf {\bibinfo {volume} {74}},\ \bibinfo {pages} {201--219} (\bibinfo {year} {2005})}\BibitemShut {NoStop}%
\bibitem [{\citenamefont {Ferracina}\ and\ \citenamefont {Spijker}(2008)}]{ferracina2008strong}%
  \BibitemOpen
  \bibfield  {author} {\bibinfo {author} {\bibfnamefont {L.}~\bibnamefont {Ferracina}}\ and\ \bibinfo {author} {\bibfnamefont {M.}~\bibnamefont {Spijker}},\ }\bibfield  {title} {\enquote {\bibinfo {title} {Strong stability of singly-diagonally-implicit runge--kutta methods},}\ }\href@noop {} {\bibfield  {journal} {\bibinfo  {journal} {Applied Numerical Mathematics}\ }\textbf {\bibinfo {volume} {58}},\ \bibinfo {pages} {1675--1686} (\bibinfo {year} {2008})}\BibitemShut {NoStop}%
\bibitem [{\citenamefont {Gottlieb}, \citenamefont {Ketcheson},\ and\ \citenamefont {Shu}(2011)}]{gottlieb2011strong}%
  \BibitemOpen
  \bibfield  {author} {\bibinfo {author} {\bibfnamefont {S.}~\bibnamefont {Gottlieb}}, \bibinfo {author} {\bibfnamefont {D.}~\bibnamefont {Ketcheson}}, \ and\ \bibinfo {author} {\bibfnamefont {C.-W.}\ \bibnamefont {Shu}},\ }\href@noop {} {\emph {\bibinfo {title} {Strong stability preserving Runge-Kutta and multistep time discretizations}}}\ (\bibinfo  {publisher} {World Scientific},\ \bibinfo {year} {2011})\BibitemShut {NoStop}%
\bibitem [{\citenamefont {Gottlieb}\ and\ \citenamefont {Shu}(1998)}]{gottlieb1998total}%
  \BibitemOpen
  \bibfield  {author} {\bibinfo {author} {\bibfnamefont {S.}~\bibnamefont {Gottlieb}}\ and\ \bibinfo {author} {\bibfnamefont {C.-W.}\ \bibnamefont {Shu}},\ }\bibfield  {title} {\enquote {\bibinfo {title} {Total variation diminishing runge-kutta schemes},}\ }\href@noop {} {\bibfield  {journal} {\bibinfo  {journal} {Mathematics of computation}\ }\textbf {\bibinfo {volume} {67}},\ \bibinfo {pages} {73--85} (\bibinfo {year} {1998})}\BibitemShut {NoStop}%
\bibitem [{\citenamefont {Gottlieb}, \citenamefont {Shu},\ and\ \citenamefont {Tadmor}(2001)}]{gottlieb2001strong}%
  \BibitemOpen
  \bibfield  {author} {\bibinfo {author} {\bibfnamefont {S.}~\bibnamefont {Gottlieb}}, \bibinfo {author} {\bibfnamefont {C.-W.}\ \bibnamefont {Shu}}, \ and\ \bibinfo {author} {\bibfnamefont {E.}~\bibnamefont {Tadmor}},\ }\bibfield  {title} {\enquote {\bibinfo {title} {Strong stability-preserving high-order time discretization methods},}\ }\href@noop {} {\bibfield  {journal} {\bibinfo  {journal} {SIAM review}\ }\textbf {\bibinfo {volume} {43}},\ \bibinfo {pages} {89--112} (\bibinfo {year} {2001})}\BibitemShut {NoStop}%
\bibitem [{\citenamefont {Higueras}(2004)}]{higueras2004strong}%
  \BibitemOpen
  \bibfield  {author} {\bibinfo {author} {\bibfnamefont {I.}~\bibnamefont {Higueras}},\ }\bibfield  {title} {\enquote {\bibinfo {title} {On strong stability preserving time discretization methods},}\ }\href@noop {} {\bibfield  {journal} {\bibinfo  {journal} {Journal of Scientific Computing}\ }\textbf {\bibinfo {volume} {21}},\ \bibinfo {pages} {193--223} (\bibinfo {year} {2004})}\BibitemShut {NoStop}%
\bibitem [{\citenamefont {Higueras}(2005)}]{higueras2005representations}%
  \BibitemOpen
  \bibfield  {author} {\bibinfo {author} {\bibfnamefont {I.}~\bibnamefont {Higueras}},\ }\bibfield  {title} {\enquote {\bibinfo {title} {Representations of runge--kutta methods and strong stability preserving methods},}\ }\href@noop {} {\bibfield  {journal} {\bibinfo  {journal} {SIAM journal on numerical analysis}\ }\textbf {\bibinfo {volume} {43}},\ \bibinfo {pages} {924--948} (\bibinfo {year} {2005})}\BibitemShut {NoStop}%
\bibitem [{\citenamefont {Hundsdorfer}, \citenamefont {Ruuth},\ and\ \citenamefont {Spiteri}(2003)}]{hundsdorfer2003monotonicity}%
  \BibitemOpen
  \bibfield  {author} {\bibinfo {author} {\bibfnamefont {W.}~\bibnamefont {Hundsdorfer}}, \bibinfo {author} {\bibfnamefont {S.~J.}\ \bibnamefont {Ruuth}}, \ and\ \bibinfo {author} {\bibfnamefont {R.~J.}\ \bibnamefont {Spiteri}},\ }\bibfield  {title} {\enquote {\bibinfo {title} {Monotonicity-preserving linear multistep methods},}\ }\href@noop {} {\bibfield  {journal} {\bibinfo  {journal} {SIAM Journal on Numerical Analysis}\ }\textbf {\bibinfo {volume} {41}},\ \bibinfo {pages} {605--623} (\bibinfo {year} {2003})}\BibitemShut {NoStop}%
\bibitem [{\citenamefont {Ketcheson}(2008)}]{ketcheson2008highly}%
  \BibitemOpen
  \bibfield  {author} {\bibinfo {author} {\bibfnamefont {D.~I.}\ \bibnamefont {Ketcheson}},\ }\bibfield  {title} {\enquote {\bibinfo {title} {Highly efficient strong stability preserving runge-kutta methods with low-storage implementations},}\ }\href@noop {} {\bibfield  {journal} {\bibinfo  {journal} {SIAM Journal on Scientific Computing}\ }\textbf {\bibinfo {volume} {30}},\ \bibinfo {pages} {2113--2136} (\bibinfo {year} {2008})}\BibitemShut {NoStop}%
\bibitem [{\citenamefont {Ketcheson}, \citenamefont {Macdonald},\ and\ \citenamefont {Gottlieb}(2009)}]{ketcheson2009optimal}%
  \BibitemOpen
  \bibfield  {author} {\bibinfo {author} {\bibfnamefont {D.~I.}\ \bibnamefont {Ketcheson}}, \bibinfo {author} {\bibfnamefont {C.~B.}\ \bibnamefont {Macdonald}}, \ and\ \bibinfo {author} {\bibfnamefont {S.}~\bibnamefont {Gottlieb}},\ }\bibfield  {title} {\enquote {\bibinfo {title} {Optimal implicit strong stability preserving runge--kutta methods},}\ }\href@noop {} {\bibfield  {journal} {\bibinfo  {journal} {Applied Numerical Mathematics}\ }\textbf {\bibinfo {volume} {59}},\ \bibinfo {pages} {373--392} (\bibinfo {year} {2009})}\BibitemShut {NoStop}%
\bibitem [{\citenamefont {Ketcheson}(2009)}]{ketcheson2009computation}%
  \BibitemOpen
  \bibfield  {author} {\bibinfo {author} {\bibfnamefont {D.}~\bibnamefont {Ketcheson}},\ }\bibfield  {title} {\enquote {\bibinfo {title} {Computation of optimal monotonicity preserving general linear methods},}\ }\href@noop {} {\bibfield  {journal} {\bibinfo  {journal} {Mathematics of Computation}\ }\textbf {\bibinfo {volume} {78}},\ \bibinfo {pages} {1497--1513} (\bibinfo {year} {2009})}\BibitemShut {NoStop}%
\bibitem [{\citenamefont {Ketcheson}(2011)}]{ketcheson2011step}%
  \BibitemOpen
  \bibfield  {author} {\bibinfo {author} {\bibfnamefont {D.~I.}\ \bibnamefont {Ketcheson}},\ }\bibfield  {title} {\enquote {\bibinfo {title} {Step sizes for strong stability preservation with downwind-biased operators},}\ }\href@noop {} {\bibfield  {journal} {\bibinfo  {journal} {SIAM Journal on Numerical Analysis}\ }\textbf {\bibinfo {volume} {49}},\ \bibinfo {pages} {1649--1660} (\bibinfo {year} {2011})}\BibitemShut {NoStop}%
\bibitem [{\citenamefont {Kubatko}, \citenamefont {Yeager},\ and\ \citenamefont {Ketcheson}(2014)}]{kubatko2014optimal}%
  \BibitemOpen
  \bibfield  {author} {\bibinfo {author} {\bibfnamefont {E.~J.}\ \bibnamefont {Kubatko}}, \bibinfo {author} {\bibfnamefont {B.~A.}\ \bibnamefont {Yeager}}, \ and\ \bibinfo {author} {\bibfnamefont {D.~I.}\ \bibnamefont {Ketcheson}},\ }\bibfield  {title} {\enquote {\bibinfo {title} {Optimal strong-stability-preserving runge--kutta time discretizations for discontinuous galerkin methods},}\ }\href@noop {} {\bibfield  {journal} {\bibinfo  {journal} {Journal of Scientific Computing}\ }\textbf {\bibinfo {volume} {60}},\ \bibinfo {pages} {313--344} (\bibinfo {year} {2014})}\BibitemShut {NoStop}%
\bibitem [{\citenamefont {Del~Buono}\ and\ \citenamefont {Lopez}(2003)}]{del2003survey}%
  \BibitemOpen
  \bibfield  {author} {\bibinfo {author} {\bibfnamefont {N.}~\bibnamefont {Del~Buono}}\ and\ \bibinfo {author} {\bibfnamefont {L.}~\bibnamefont {Lopez}},\ }\bibfield  {title} {\enquote {\bibinfo {title} {A survey on methods for computing matrix exponentials in numerical schemes for odes},}\ }in\ \href@noop {} {\emph {\bibinfo {booktitle} {International Conference on Computational Science}}}\ (\bibinfo {organization} {Springer},\ \bibinfo {year} {2003})\ pp.\ \bibinfo {pages} {111--120}\BibitemShut {NoStop}%
\bibitem [{\citenamefont {Ruiz}\ \emph {et~al.}(2016)\citenamefont {Ruiz}, \citenamefont {Sastre}, \citenamefont {Ibáñez},\ and\ \citenamefont {Defez}}]{RUIZ2016370}%
  \BibitemOpen
  \bibfield  {author} {\bibinfo {author} {\bibfnamefont {P.}~\bibnamefont {Ruiz}}, \bibinfo {author} {\bibfnamefont {J.}~\bibnamefont {Sastre}}, \bibinfo {author} {\bibfnamefont {J.}~\bibnamefont {Ibáñez}}, \ and\ \bibinfo {author} {\bibfnamefont {E.}~\bibnamefont {Defez}},\ }\bibfield  {title} {\enquote {\bibinfo {title} {High performance computing of the matrix exponential},}\ }\href {\doibase https://doi.org/10.1016/j.cam.2015.04.001} {\bibfield  {journal} {\bibinfo  {journal} {Journal of Computational and Applied Mathematics}\ }\textbf {\bibinfo {volume} {291}},\ \bibinfo {pages} {370--379} (\bibinfo {year} {2016})},\ \bibinfo {note} {mathematical Modeling and Computational Methods}\BibitemShut {NoStop}%
\bibitem [{\citenamefont {Arioli}, \citenamefont {Codenotti},\ and\ \citenamefont {Fassino}(1996)}]{ARIOLI1996111}%
  \BibitemOpen
  \bibfield  {author} {\bibinfo {author} {\bibfnamefont {M.}~\bibnamefont {Arioli}}, \bibinfo {author} {\bibfnamefont {B.}~\bibnamefont {Codenotti}}, \ and\ \bibinfo {author} {\bibfnamefont {C.}~\bibnamefont {Fassino}},\ }\bibfield  {title} {\enquote {\bibinfo {title} {The padé method for computing the matrix exponential},}\ }\href {\doibase https://doi.org/10.1016/0024-3795(94)00190-1} {\bibfield  {journal} {\bibinfo  {journal} {Linear Algebra and its Applications}\ }\textbf {\bibinfo {volume} {240}},\ \bibinfo {pages} {111--130} (\bibinfo {year} {1996})}\BibitemShut {NoStop}%
\bibitem [{\citenamefont {van~der Pol}(1926)}]{VanderPolModel}%
  \BibitemOpen
  \bibfield  {author} {\bibinfo {author} {\bibfnamefont {B.}~\bibnamefont {van~der Pol}},\ }\bibfield  {title} {\enquote {\bibinfo {title} {On relaxation-oscillations},}\ }\href {\doibase 10.1080/14786442608564127} {\bibfield  {journal} {\bibinfo  {journal} {The London, Edinburgh, and Dublin Philosophical Magazine and Journal of Science}\ }\textbf {\bibinfo {volume} {2}},\ \bibinfo {pages} {978--992} (\bibinfo {year} {1926})}\BibitemShut {NoStop}%
\bibitem [{\citenamefont {Schäfer}(1975)}]{Schaefer1975HIR}%
  \BibitemOpen
  \bibfield  {author} {\bibinfo {author} {\bibfnamefont {E.}~\bibnamefont {Schäfer}},\ }\bibfield  {title} {\enquote {\bibinfo {title} {A new approach to explain the “high irradiance responses” of photomorphogenesis on the basis of phytochrome},}\ }\href {\doibase 10.1007/BF00276015} {\bibfield  {journal} {\bibinfo  {journal} {Journal of Mathematical Biology}\ }\textbf {\bibinfo {volume} {2}},\ \bibinfo {pages} {41--56} (\bibinfo {year} {1975})}\BibitemShut {NoStop}%
\bibitem [{\citenamefont {Robertson}(1966)}]{robertson1966reaction}%
  \BibitemOpen
  \bibfield  {author} {\bibinfo {author} {\bibfnamefont {H.}~\bibnamefont {Robertson}},\ }\bibfield  {title} {\enquote {\bibinfo {title} {The solution of a set of reaction rate equations},}\ }in\ \href@noop {} {\emph {\bibinfo {booktitle} {Numerical Analysis: An Introduction}}}\ (\bibinfo  {publisher} {Academic Press},\ \bibinfo {year} {1966})\ pp.\ \bibinfo {pages} {178--182}\BibitemShut {NoStop}%
\bibitem [{\citenamefont {Virtanen}\ \emph {et~al.}(2020)\citenamefont {Virtanen}, \citenamefont {Gommers}, \citenamefont {Oliphant}, \citenamefont {Haberland}, \citenamefont {Reddy}, \citenamefont {Cournapeau}, \citenamefont {Burovski}, \citenamefont {Peterson}, \citenamefont {Weckesser}, \citenamefont {Bright}, \citenamefont {{van der Walt}}, \citenamefont {Brett}, \citenamefont {Wilson}, \citenamefont {Millman}, \citenamefont {Mayorov}, \citenamefont {Nelson}, \citenamefont {Jones}, \citenamefont {Kern}, \citenamefont {Larson}, \citenamefont {Carey}, \citenamefont {Polat}, \citenamefont {Feng}, \citenamefont {Moore}, \citenamefont {{VanderPlas}}, \citenamefont {Laxalde}, \citenamefont {Perktold}, \citenamefont {Cimrman}, \citenamefont {Henriksen}, \citenamefont {Quintero}, \citenamefont {Harris}, \citenamefont {Archibald}, \citenamefont {Ribeiro}, \citenamefont {Pedregosa}, \citenamefont {{van Mulbregt}},\ and\ \citenamefont {{SciPy 1.0 Contributors}}}]{2020SciPy-NMeth}%
  \BibitemOpen
  \bibfield  {author} {\bibinfo {author} {\bibfnamefont {P.}~\bibnamefont {Virtanen}}, \bibinfo {author} {\bibfnamefont {R.}~\bibnamefont {Gommers}}, \bibinfo {author} {\bibfnamefont {T.~E.}\ \bibnamefont {Oliphant}}, \bibinfo {author} {\bibfnamefont {M.}~\bibnamefont {Haberland}}, \bibinfo {author} {\bibfnamefont {T.}~\bibnamefont {Reddy}}, \bibinfo {author} {\bibfnamefont {D.}~\bibnamefont {Cournapeau}}, \bibinfo {author} {\bibfnamefont {E.}~\bibnamefont {Burovski}}, \bibinfo {author} {\bibfnamefont {P.}~\bibnamefont {Peterson}}, \bibinfo {author} {\bibfnamefont {W.}~\bibnamefont {Weckesser}}, \bibinfo {author} {\bibfnamefont {J.}~\bibnamefont {Bright}}, \bibinfo {author} {\bibfnamefont {S.~J.}\ \bibnamefont {{van der Walt}}}, \bibinfo {author} {\bibfnamefont {M.}~\bibnamefont {Brett}}, \bibinfo {author} {\bibfnamefont {J.}~\bibnamefont {Wilson}}, \bibinfo {author} {\bibfnamefont {K.~J.}\ \bibnamefont {Millman}}, \bibinfo {author} {\bibfnamefont {N.}~\bibnamefont {Mayorov}}, \bibinfo {author} {\bibfnamefont
  {A.~R.~J.}\ \bibnamefont {Nelson}}, \bibinfo {author} {\bibfnamefont {E.}~\bibnamefont {Jones}}, \bibinfo {author} {\bibfnamefont {R.}~\bibnamefont {Kern}}, \bibinfo {author} {\bibfnamefont {E.}~\bibnamefont {Larson}}, \bibinfo {author} {\bibfnamefont {C.~J.}\ \bibnamefont {Carey}}, \bibinfo {author} {\bibfnamefont {{\.I}.}~\bibnamefont {Polat}}, \bibinfo {author} {\bibfnamefont {Y.}~\bibnamefont {Feng}}, \bibinfo {author} {\bibfnamefont {E.~W.}\ \bibnamefont {Moore}}, \bibinfo {author} {\bibfnamefont {J.}~\bibnamefont {{VanderPlas}}}, \bibinfo {author} {\bibfnamefont {D.}~\bibnamefont {Laxalde}}, \bibinfo {author} {\bibfnamefont {J.}~\bibnamefont {Perktold}}, \bibinfo {author} {\bibfnamefont {R.}~\bibnamefont {Cimrman}}, \bibinfo {author} {\bibfnamefont {I.}~\bibnamefont {Henriksen}}, \bibinfo {author} {\bibfnamefont {E.~A.}\ \bibnamefont {Quintero}}, \bibinfo {author} {\bibfnamefont {C.~R.}\ \bibnamefont {Harris}}, \bibinfo {author} {\bibfnamefont {A.~M.}\ \bibnamefont {Archibald}}, \bibinfo {author}
  {\bibfnamefont {A.~H.}\ \bibnamefont {Ribeiro}}, \bibinfo {author} {\bibfnamefont {F.}~\bibnamefont {Pedregosa}}, \bibinfo {author} {\bibfnamefont {P.}~\bibnamefont {{van Mulbregt}}}, \ and\ \bibinfo {author} {\bibnamefont {{SciPy 1.0 Contributors}}},\ }\bibfield  {title} {\enquote {\bibinfo {title} {{{SciPy} 1.0: Fundamental Algorithms for Scientific Computing in Python}},}\ }\href {\doibase 10.1038/s41592-019-0686-2} {\bibfield  {journal} {\bibinfo  {journal} {Nature Methods}\ }\textbf {\bibinfo {volume} {17}},\ \bibinfo {pages} {261--272} (\bibinfo {year} {2020})}\BibitemShut {NoStop}%
\bibitem [{\citenamefont {FitzHugh}(1961)}]{FITZHUGH1961445}%
  \BibitemOpen
  \bibfield  {author} {\bibinfo {author} {\bibfnamefont {R.}~\bibnamefont {FitzHugh}},\ }\bibfield  {title} {\enquote {\bibinfo {title} {Impulses and physiological states in theoretical models of nerve membrane},}\ }\href {\doibase https://doi.org/10.1016/S0006-3495(61)86902-6} {\bibfield  {journal} {\bibinfo  {journal} {Biophysical Journal}\ }\textbf {\bibinfo {volume} {1}},\ \bibinfo {pages} {445--466} (\bibinfo {year} {1961})}\BibitemShut {NoStop}%
\bibitem [{\citenamefont {Nagumo}, \citenamefont {Arimoto},\ and\ \citenamefont {Yoshizawa}(1962)}]{4066548}%
  \BibitemOpen
  \bibfield  {author} {\bibinfo {author} {\bibfnamefont {J.}~\bibnamefont {Nagumo}}, \bibinfo {author} {\bibfnamefont {S.}~\bibnamefont {Arimoto}}, \ and\ \bibinfo {author} {\bibfnamefont {S.}~\bibnamefont {Yoshizawa}},\ }\bibfield  {title} {\enquote {\bibinfo {title} {An active pulse transmission line simulating nerve axon},}\ }\href {\doibase 10.1109/JRPROC.1962.288235} {\bibfield  {journal} {\bibinfo  {journal} {Proceedings of the IRE}\ }\textbf {\bibinfo {volume} {50}},\ \bibinfo {pages} {2061--2070} (\bibinfo {year} {1962})}\BibitemShut {NoStop}%
\bibitem [{\citenamefont {Cartwright}\ \emph {et~al.}(1999)\citenamefont {Cartwright}, \citenamefont {EGU\'{I}LUZ}, \citenamefont {HERN\'{A}NDEZ-GARC\'{I}A},\ and\ \citenamefont {PIRO}}]{doi:10.1142/S0218127499001620}%
  \BibitemOpen
  \bibfield  {author} {\bibinfo {author} {\bibfnamefont {J.~H.~E.}\ \bibnamefont {Cartwright}}, \bibinfo {author} {\bibfnamefont {V.~M.}\ \bibnamefont {EGU\'{I}LUZ}}, \bibinfo {author} {\bibfnamefont {E.}~\bibnamefont {HERN\'{A}NDEZ-GARC\'{I}A}}, \ and\ \bibinfo {author} {\bibfnamefont {O.}~\bibnamefont {PIRO}},\ }\bibfield  {title} {\enquote {\bibinfo {title} {Dynamics of elastic excitable media},}\ }\href {\doibase 10.1142/S0218127499001620} {\bibfield  {journal} {\bibinfo  {journal} {International Journal of Bifurcation and Chaos}\ }\textbf {\bibinfo {volume} {09}},\ \bibinfo {pages} {2197--2202} (\bibinfo {year} {1999})}\BibitemShut {NoStop}%
\bibitem [{\citenamefont {Lucero}\ and\ \citenamefont {Schoentgen}(2013)}]{doi:10.1121/1.4798467}%
  \BibitemOpen
  \bibfield  {author} {\bibinfo {author} {\bibfnamefont {J.~C.}\ \bibnamefont {Lucero}}\ and\ \bibinfo {author} {\bibfnamefont {J.}~\bibnamefont {Schoentgen}},\ }\bibfield  {title} {\enquote {\bibinfo {title} {Modeling vocal fold asymmetries with coupled van der pol oscillators},}\ }\href {\doibase 10.1121/1.4798467} {\bibfield  {journal} {\bibinfo  {journal} {Proceedings of Meetings on Acoustics}\ }\textbf {\bibinfo {volume} {19}},\ \bibinfo {pages} {060165} (\bibinfo {year} {2013})}\BibitemShut {NoStop}%
\end{thebibliography}%

\end{document}